\documentclass[reqno,oneside]{amsart} 

\usepackage{amsthm}

\usepackage{
amssymb,amsmath,enumerate,stackrel,
mathrsfs,verbatim,stmaryrd,
leftidx}
\usepackage{xcolor} 
\usepackage{xr-hyper}
\usepackage{nicefrac,stackrel}
\usepackage[all,cmtip]{xy}
\usepackage[utf8]{inputenc} 
\chardef\cprime"7E 
\usepackage{cases}

\usepackage[T1]{fontenc} 

\usepackage[pagebackref,hyperindex,linktocpage=true]{hyperref}
\hypersetup{
    colorlinks,
    linkcolor={red!50!black},
    citecolor={red!50!black}, 
    urlcolor={red!50!black}, 
    filecolor={red!50!black} 
}

\usepackage{tikz}
\usetikzlibrary{matrix,arrows}

\usepackage{xcolor}

\usepackage{theoremref}
\definecolor{labelkey}{rgb}{1,0,0}

\makeatletter
\@addtoreset{equation}{section}
\makeatother

\numberwithin{equation}{section}

\theoremstyle{definition}
\newtheorem{Defi}{Definition}[section] \newcommand{\defi}{\begin{Defi}} \newcommand{\xdefi}{\end{Defi}} 
\newtheorem{DefiLemm}[Defi]{Definition and Lemma} \newcommand{\defilemm}{\begin{DefiLemm}} \newcommand{\xdefilemm}{\end{DefiLemm}} 
\newtheorem{Bsp}[Defi]{Example} \newcommand{\exam}{\begin{Bsp}} \newcommand{\xexam}{\end{Bsp}} 
\newtheorem{Syno}[Defi]{Synopsis} \newcommand{\syno}{\begin{Syno}} \newcommand{\xsyno}{\end{Syno}} 
\newtheorem{Bem}[Defi]{Remark} \newcommand{\rema}{\begin{Bem}} \newcommand{\xrema}{\end{Bem}} 
\newtheorem{Notation}[Defi]{Notation} \newcommand{\nota}{\begin{Notation}} \newcommand{\xnota}{\end{Notation}} 
\newtheorem{Warning}[Defi]{Warning} \newcommand{\warn}{\begin{Warning}} \newcommand{\xwarn}{\end{Warning}} 
\newtheorem{Situation}[Defi]{Situation} \newcommand{\situ}{\begin{Situation}} \newcommand{\xsitu}{\end{Situation}} 

\theoremstyle{plain}
\newtheorem{Theo}[Defi]{Theorem} \newcommand{\theo}{\begin{Theo}} \newcommand{\xtheo}{\end{Theo}} 
\newtheorem{Satz}[Defi]{Proposition} \newcommand{\prop}{\begin{Satz}} \newcommand{\xprop}{\end{Satz}} 
\newtheorem{Lemm}[Defi]{Lemma} \newcommand{\lemm}{\begin{Lemm}} \newcommand{\xlemm}{\end{Lemm}} 
\newtheorem{Coro}[Defi]{Corollary} \newcommand{\coro}{\begin{Coro}} \newcommand{\xcoro}{\end{Coro}}
\newtheorem{Ques}[Defi]{Question} \newcommand{\ques}{\begin{Ques}} \newcommand{\xques}{\end{Ques}}
\newtheorem{Conj}[Defi]{Conjecture} \newcommand{\conj}{\begin{Conj}} \newcommand{\xconj}{\end{Conj}}

\newcommand{\refsect}[1]{Section \ref{sect--#1}}
\newcommand{\refit}[1]{(\ref{item--#1})}
\newcommand{\refeq}[1]{(\ref{eqn--#1})}

\newcommand{\eqn}{\begin{equation}} \newcommand{\xeqn}{\end{equation}}
\newcommand{\eqnarr}{\begin{eqnarray*}} \newcommand{\xeqnarr}{\end{eqnarray*}}
\newcommand{\eqnarra}{\begin{eqnarray}} \newcommand{\xeqnarra}{\end{eqnarray}}

\newcommand{\pf}{\begin{proof}} \newcommand{\xpf}{\end{proof}}

\usepackage[top=1.5cm, bottom=1.5cm, left=2cm, right=2cm, heightrounded, marginparwidth=3.5cm, marginparsep=.2cm]{geometry}

\usepackage{color}

\newcommand{\nc}{\newcommand}
\nc{\StP}[1]{\cite[Tag~\href{http://stacks.math.columbia.edu/tag/#1}{#1}]{StacksProject}}
\nc{\StPd}[2]{\cite[Tags~\href{http://stacks.math.columbia.edu/tag/#1}{#1}, \href{http://stacks.math.columbia.edu/tag/#2}{#2}]{StacksProject}} 

\nc{\on}{\operatorname}
\nc{\aff}{{\on{aff}}}
\nc{\modi}{{\on{mod}}} 
\nc{\even}{{\on{even}}}
\nc{\odd}{{\on{odd}}}
\nc{\naive}{{\on{naive}}}
\nc{\hofib}{\on{hofib}}
\nc{\Bun}{\on{Bun}}
\nc{\ad}{{\on{ad}}}
\nc{\lft}{{\on{lft}}}

\nc{\Weil}{{\on{Weil}}} 
\nc{\FWeil}{{\on{FWeil}}} 
\nc{\cons}{{\on{cons}}} 
\nc{\tot}{{\on{Tot}}} 

\nc{\str}{\on{-}}
\nc{\perf}{{\on{perf}}}
\nc{\Rel}{{\on{Pos}}}
\nc{\lan}{\langle}
\nc{\ran}{\rangle}

\nc{\bbA}{{\mathbb A}} 
\nc{\bbB}{{\mathbb B}}
\nc{\bbC}{{\mathbb C}}
\nc{\bbD}{{\mathbb D}}
\nc{\bbE}{{\mathbb E}}
\nc{\bbF}{{\mathbb F}}
\nc{\bbG}{{\mathbf G}}
\nc{\bbH}{{\mathbb H}}
\nc{\bbI}{{\mathbb I}}
\nc{\bbJ}{{\mathbb J}}
\nc{\bbK}{{\mathbb K}}
\nc{\bbL}{{\mathbb L}}
\nc{\bbM}{{\mathbb M}}
\nc{\bbN}{{\N}} 
\nc{\bbO}{{\mathbb O}}
\nc{\bbP}{{\mathbb P}} 
\nc{\bbQ}{{\mathbb Q}} 
\nc{\bbR}{{\mathbb R}}
\nc{\bbS}{{\mathbb S}}
\nc{\bbT}{{\mathbb T}}
\nc{\bbU}{{\mathbb U}}
\nc{\bbV}{{\mathbb V}}
\nc{\bbW}{{\mathbb W}}
\nc{\bbX}{{\mathbb X}}
\nc{\bbY}{{\mathbb Y}}
\nc{\bbZ}{{\mathbb Z}}

\nc{\calA}{{\mathcal A}}
\nc{\calB}{{\mathcal B}}
\nc{\calC}{{\mathcal C}}
\nc{\calD}{{\mathcal D}}
\nc{\calE}{{\mathcal E}}
\nc{\calF}{{\mathcal F}}
\nc{\calG}{{\mathcal G}}
\nc{\calH}{{\mathcal H}}
\nc{\calI}{{\mathcal I}}
\nc{\calJ}{{\mathcal J}}
\nc{\calK}{{\mathcal K}}
\nc{\calL}{{\mathcal L}}
\nc{\calM}{{\mathcal M}}
\nc{\calN}{{\mathcal N}}
\nc{\calO}{{\mathcal O}}
\nc{\calP}{{\mathcal P}}
\nc{\calQ}{{\mathcal Q}}
\nc{\calR}{{\mathcal R}}
\nc{\calS}{{\mathcal S}}
\nc{\calT}{{\mathcal T}}
\nc{\calU}{{\mathcal U}}
\nc{\calV}{{\mathcal V}}
\nc{\calW}{{\mathcal W}}
\nc{\calX}{{\mathcal X}}
\nc{\calY}{{\mathcal Y}}
\nc{\calZ}{{\mathcal Z}}

\nc{\Sht}{{\on{Sht}}}
\nc{\Frob}{{\on{Frob}}}
\nc{\Hecke}{{\on{Hecke}}}
\nc{\inv}{{\on{inv}}}
\nc{\Conv}{{\on{Conv}}}
\nc{\triv}{{\on{triv}}}
\nc{\Isom}{{\on{Isom}}}

\nc{\scrB}{{\mathscr{B}}}
\nc{\scrA}{{\mathscr{A}}}
\nc{\bbf}{{\mathbf{f}}}
\nc{\bba}{{\mathbf{a}}}
\nc{\rig}{{\mathrm rig}}

\nc{\Indft}{{\pi_1\!\on{-}\!\mathrm{Indft}}}
\nc{\IndPerf}{{\pi_1\!\on{-}\!\mathrm{IndPerf}}}
\nc{\AS}{{\on{AS}}}

\nc{\al}{\alpha}
\nc{\be}{\beta}
\nc{\ga}{\gamma}
\nc{\la}{\lambda}
\nc{\qcqs}{{\on{qcqs}}}
\nc{\Bmu}{{\boldsymbol \mu}}

\nc{\pot}[1]{ [\hspace{-0,5mm}[ {#1} ]\hspace{-0,5mm}] }
\nc{\rpot}[1]{ (\hspace{-0,7mm}( {#1} )\hspace{-0,7mm}) }

\nc{\defined}{\hspace{0.1cm}\stackrel{\text{\tiny \rm def}}{=}\hspace{0.1cm}}
\nc{\co}{\colon}
\nc{\specto}{{\leadsto}}

\newcommand{\category}[1]{\mathrm{#1}}
\newcommand{\DGCat}{\category{DGCat}} 
\newcommand{\cont}{\category{cont}} 
\newcommand{\Cat}{\category{Cat}} 
\newcommand{\Ex}{\category{Ex}} 
\newcommand{\CatEx}{\Cat^{\Ex}} 
\newcommand{\Fun}{\category{Fun}} 
\newcommand{\Ani}{\category{Ani}} 
\newcommand{\cts}{\mathrm{cts}}  
\newcommand{\Zar}{\mathrm{Zar}} 
\renewcommand{\Pr}{\category{Pr}}
\newcommand{\PrL}{\Pr^\category{L}} 
\newcommand{\PrSt}{\Pr^{\category{St}}} 
\newcommand{\PrStL}{\PrSt_\Lambda} 
\newcommand{\Ind}{\category{Ind}} 
\newcommand{\Mod}{\category{Mod}} 
\newcommand{\Sp}{\category{Sp}} 
\newcommand{\Perf}{\category{Perf}} 
\newcommand{\Picond}{\Pi_{\mathrm{cond}}} 
\newcommand{\Tot}{\mathrm {Tot}} 
\newcommand{\Idem}{\mathrm {Idem}} 
\newcommand{\Dbc}{\mathrm{D}_{\mathrm{c}}^{\mathrm{b}}}
\newcommand{\Dcons}{\mathrm{D}_\cons}
\newcommand{\Dlis}{\mathrm{D}_\lis}

\def\fp{\mathrm{fp}} 

\font\tencyr=wncyr10
\font\sevencyr=wncyr7
\font\fivecyr=wncyr5
\newcommand{\colim}{\operatornamewithlimits{colim}} 
\def\id{{\rm id}} 
\def\To#1#2{\mathop{\count0=#1 \loop\ifnum\count0>0 \smash-\mkern-7mu \advance\count0 -1 \repeat \mathord\rightarrow}\limits^{#2}} 
\def\Maps{\mathop{\rm Maps}\nolimits} 
\def\Hom{\mathop{\rm Hom}\nolimits} 

\def\Ind{\category{Ind}} 
\def\RHom{\mathop{\rm RHom}\nolimits} 
\def\Sht{\mathop{\rm Sht}\nolimits} 
\def\IHom{\underline{\Hom}} 
\def\Aut{\mathop{\rm Aut}\nolimits} 
\def\et{\mathrm{\acute et}} 
\def\proet{\mathrm{pro\acute et}} 

\definecolor{hellgrau}{RGB}{200,200,200} 
\definecolor{dunkelgrau}{RGB}{160,160,160} 
\definecolor{hellblau}{RGB}{194, 215, 249} %
\definecolor{dunkelblau}{RGB}{68, 128, 226} %

\def\Z{{\bbZ}} 
\def\Fp{{{\mathbb F}_p}} %
\def\N{{\mathbb N}} 
\def\Q{{\bbQ}} 
\def\Ql{{\Q_\ell}} 
\def\Zl{{\Z_\ell}} 

\def\H{{\rm H}} 
\def\ii{$\infty$}

\def\Gal{{\rm Gal}} 
\def\ctf{{\rm ctf}} 
\def\Cond{{\rm Cond}} 

\def\lis{{\rm lis}} 
\def\indlis{{\rm indlis}} 
\def\indcons{{\rm indcons}} 

\def\RG{\R \Gamma} 

\def\Spec{\mathop{\rm Spec}} 
\newcommand{\D}{\category{D}} 

\def\R{{\rm R}} 
\def\sbuildrel#1\over#2{\mathrel{\smash{\mathop{\kern0pt #2}\limits^{#1}}}}

\let\x\times

\renewcommand{\t}{\otimes}

\renewcommand{\r}{\rightarrow}
\newcommand{\lr}{\longrightarrow}
\newcommand{\hr}{\hookrightarrow}

\def\matrix#1{\null\,\vcenter{\normalbaselines
    \ialign{\hfil$##$\hfil&&\quad\hfil$##$\hfil\crcr
      \mathstrut\crcr\noalign{\kern-\baselineskip}
      #1\crcr\mathstrut\crcr\noalign{\kern-\baselineskip}}}\,}

\newdimen\harrowsize
\def\mapright#1{\smash{\mathop{\hbox to\harrowsize{\rightarrowfill}}\limits^{#1}}}

{\catcode`@=11
\gdef\cal{\fam\tw@}
\global\let\over\@@over
\global\let\atop\@@atop
\global\let\above\@@above
\global\let\overwithdelims\@@overwithdelims
\global\let\atopwithdelims\@@atopwithdelims
\global\let\abovewithdelims\@@abovewithdelims
\gdef\eqalign#1{\null\,\vcenter{\openup\jot\m@th
\ialign{\strut\hfil$\displaystyle{##}$&$\displaystyle{{}##}$\hfil
      \crcr#1\crcr}}\,}
\newskip\xcentering \global\xcentering=0pt plus 1000pt minus 1000pt
\gdef\eqalignno#1{\displ@y \tabskip\xcentering
  \halign to\displaywidth{\hfil$\@lign\displaystyle{##}$\tabskip\z@skip
    &$\@lign\displaystyle{{}##}$\hfil\tabskip\xcentering
    &\llap{$\@lign##$}\tabskip\z@skip\crcr
    #1\crcr}}
\global\def\cases#1{\left\{\,\vcenter{\normalbaselines\m@th
    \ialign{$##\hfil$&\quad##\hfil\crcr#1\crcr}}\right.}
\gdef\eqlabel#1{\refstepcounter{equation}\label{eqn--#1}\eqno\hbox{\@eqnnum}}
}

\def \journal#1
{ 
\noindent\colorbox{dunkelblau}{\parbox{\dimexpr\textwidth-2\fboxsep\relax}{#1}}
}

\begin{document}

\author{Tamir Hemo, Timo Richarz and Jakob Scholbach}

\title{Constructible sheaves on schemes}

\thanks{*The second named author T.R.~is funded by the European Research Council (ERC) under Horizon Europe (grant agreement nº 101040935), by the Deutsche Forschungsgemeinschaft (DFG, German Research Foundation) TRR 326 \textit{Geometry and Arithmetic of Uniformized Structures}, project number 444845124 and the LOEWE professorship in Algebra.
J.S.~was supported by Deutsche Forschungsgemeinschaft (DFG), EXC 2044–390685587, Mathematik Münster: Dynamik–Geometrie–Struktur. }

\address{Caltech, Department of Mathematics, 91125 Pasadena, CA, USA}
\email{themo@caltech.edu}

\address{Technical University of Darmstadt, Department of Mathematics, 64289 Darmstadt, Germany}
\email{richarz@mathematik.tu-darmstadt.de}

\address{Università degli Studi di Padova, Dipartimento di Matematica, 35139 Padova, Italia}

\begin{abstract}
We present a uniform theory of constructible sheaves on arbitrary schemes with coefficients in topological or even condensed rings.
This is accomplished by defining lisse sheaves to be the dualizable objects in the derived \ii-category of proétale sheaves, while constructible sheaves are those that are lisse on a stratification.

We show that constructible sheaves satisfy proétale descent.
We also establish a t-structure on constructible sheaves in a wide range of cases.
We finally provide a toolset to manipulate categories of constructible sheaves with respect to the choices of coefficient rings, and use this to prove that our notions reproduce and extend the various approaches to, say, constructible $\ell$-adic sheaves in the literature.
\end{abstract}

\maketitle

\tableofcontents

\section{Introduction}

Initially, étale cohomology of a scheme $X$ only works well for sheaves with torsion coefficients, such as $\Z/n$. 
In order to have a cohomology theory producing vector spaces over characteristic zero fields, one fixes a prime number $\ell$ and uses $\ell$-adic sheaves \cite[Exposé VI]{SGA5} whose cohomology groups are $\bbQ_\ell$-vector spaces.
The name ``sheaves'' is fully justified with the advent of the pro\'etale topology introduced in \cite{BhattScholze:ProEtale}: such $\ell$-adic sheaves are ordinary topos-theoretic sheaves on the site $X_\proet$.
This is decidedly more conceptual than the classical approach which builds $\ell$-adic sheaves by retracing the formula $\bbQ_\ell=(\lim_n \bbZ/\ell^n)[\ell^{-1}]$ on the level of sheaves.

For many purposes, it is useful to impose a finiteness condition on sheaves known as constructibility.
This is accomplished in \cite{BhattScholze:ProEtale} for adic coefficient rings like $\Zl = \lim \Z/\ell^n$.
Their category satisfies pro\'etale descent and compares well to \cite{SGA5, Deligne:Weil2, Ekedahl:Adic}.
If the underlying topological space of $X$ is Noetherian, \cite{BhattScholze:ProEtale} also defines a category of constructible $\Ql$-sheaves.
However, the finiteness condition on $X$ prevents one from formulating descent on $X_\proet$: typical pro\'etale covers $\{U_i\r X\}$ have a profinite set of connected components and so are not topologically Noetherian.
Examples show that the notion of pro\'etale locally constant sheaves is not well behaved in general. 
This prevents an obvious generalization of constructible $\bbQ_\ell$-sheaves to such large schemes. 

Here we introduce a notion of lisse and constructible sheaves on arbitrary schemes $X$ with coefficients in an arbitrary condensed (unital, commutative) ring $\Lambda$.
Recall from \cite{ClausenScholze:CondensedMathematics} (see also \cite{BarwickHaine:Pyknotic}) that a condensed ring is a sheaf of rings on the site $*_\proet$ of profinite sets. 
Examples abound, since every T1-topological ring naturally gives rise to a condensed ring via the Yoneda embedding.
Let $\Lambda_*=\Gamma(*,\Lambda)$ be the underlying ring.

The pullback of $\Lambda$ along the canonical map of sites $X_\proet\r *_\proet$ is a sheaf of rings on $X_\proet$. 
We denote by $\D(X,\Lambda)$ the derived \ii-category of the abelian category of $\Lambda$-sheaves on $X_\proet$.
This is a $\Lambda_*$-linear symmetric monoidal closed stable \ii-category (so its homotopy category is a triangulated category). 
Recall the notion of dualizable objects in symmetric monoidal categories \cite[Section 4.6.1]{Lurie:HA}.
For example, for any ring $R$, the dualizable objects in $\Mod_R$, the derived \ii-category of $R$-modules, are the perfect complexes, that is, bounded complexes of finite projective $R$-modules.
This subcategory is denoted by $\Perf_R \subset \Mod_R$.

\defi
\thlabel{definition.lisse.constructible.introduction}
Let $X$ be a scheme, and let $\Lambda$ be a condensed ring. 
\begin{enumerate}
\item 
A sheaf $M\in \D(X,\Lambda)$ is called {\it lisse} if it is dualizable. 

\item 
A sheaf $M\in \D(X,\Lambda)$ is called {\it constructible} if for every open affine subscheme $U\subset X$ there exists a finite subdivision $U_i\subset U$ into constructible locally closed subschemes such that each $M|_{U_i}$ is lisse. 
\end{enumerate}
\xdefi

We denote by $\D_\lis(X,\Lambda)\subset \D_\cons(X,\Lambda)$ the full subcategories of $\D(X,\Lambda)$ of lisse, respectively constructible sheaves.
By the setup, these are naturally $\Lambda_*$-linear symmetric monoidal stable \ii-categories. 

The idea that the lisse-ness condition on a sheaf can be expressed by means of dualizability is known, for example, in the context of motives or étale torsion sheaves \cite{CisinskiDeglise:Etale}, but seems to be new for coefficients such as $\Lambda = \Ql$ or $\Zl$.
The present paper offers the insight that this highly conceptual finiteness notion is also computationally manageable in the pro\'etale topology. 
The following lemma is the key stepping stone making this possible. 
Recall that every scheme admits a pro\'etale cover $\{U_i\r X\}$ where the $U_i$ are w-contractible affine schemes.
These objects form a basis of the topology on $X_\proet$ \cite[Theorem~1.5]{BhattScholze:ProEtale}. 

\lemm[\thref{dualizable.objects.contractible.lem}]
\thlabel{dualizable.objects.contractible.lem.intro}
Assume that $X$ is w-contractible affine. 
Then the global sections functor induces an equivalence of \ii-categories
\[
\RG(X,\str)\co \D_\lis(X,\Lambda)\stackrel\cong \lr \Perf_{\Gamma(X,\Lambda)}.
\]
\xlemm

Even for $X=*$ and $\Lambda=\bbZ$, this is noteworthy: the derived \ii-category of condensed abelian groups $\D(*,\bbZ)$ contains fully faithfully the category of compactly generated T1-topological abelian groups, and is thus much larger than the \ii-category $\D_\lis(*, \Z) \cong \Perf_\Z$, which consists of bounded complexes of finitely generated free abelian groups equipped with the discrete topology.

\lemm[\thref{lisse.constructible.hyperdescent.coro}]
The functors $U\mapsto \D_\lis(U,\Lambda), \D_\cons(U,\Lambda)$ are hypersheaves of \ii-categories on $X_\proet$.
\xlemm

Together with \thref{dualizable.objects.contractible.lem.intro}, this says that $U\mapsto \D_\lis(U,\Lambda)$ is the unique hypersheaf on $X_\proet$ whose values at the basis of w-contractible affines compute the \ii-category $\Perf_{\Gamma(U,\Lambda)}$.
It can be used to relate lisse sheaves to the condensed shape (\thref{lisse.sheaves.via.condensed.shape}), also known as the proétale homotopy type, which is related to the stratified shape developed in \cite{Barwick.Glasman.Haine.exodromy}.
The above functors are even sheaves in the arc topology \cite{BhattMathew:Arc} for favorable choices of $\Lambda$ such as finite discrete rings, algebraic extensions $E\supset \bbQ_\ell$ and their rings of integers $\calO_E$, see \cite[Theorem 2.2]{HansenScholze:RelativePerversity}.

By the above definition, lisse and constructible sheaves start out life in a derived setting.
The natural t-structure on $\D(X,\Lambda)$ restricts to a t-structure on such sheaves only under additional assumptions on $X$ and $\Lambda$: 

\defi 
\thlabel{t-admissible}
A condensed ring $\Lambda$ is called {\it t-admissible} if $\Lambda_*$ is regular coherent (that is, every finitely generated ideal is finitely presented and has finite projective dimension) and, for any extremally disconnected profinite set $S$, the map $\Lambda_*\r \Gamma(S,\Lambda)$ is flat. 
\xdefi

In \refsect{t.structures.sect}, we show that $\Lambda$ is t-admissible if and only if the t-structure on $\D(*,\Lambda)$ restricts a t-structure on $\D_\lis(*,\Lambda)$.
Examples of t-admissible condensed rings include discrete rings that are regular Noetherian of finite Krull dimension, and all T1-topological rings such that $\Lambda_*$ is semi-hereditary (=every finitely generated ideal is projective).
This covers algebraic field extensions $E\supset \bbQ_\ell$ and their rings of integers $\calO_E$, but also the real and complex numbers $\bbR$, $\bbC$ with their Euclidean topology and the ring of adeles $\bbA_K^T$ prime to some finite set of places $T$ in some number field $K$. 

\theo[\thref{t.structure.condensed}]
Let $\Lambda$ be a t-admissible condensed ring. 
\begin{itemize}
  \item 
If $X$ has Zariski-locally finitely many irreducible components, then
the natural t-structure on $\D(X,\Lambda)$ restricts to a t-structure on $\D_\lis(X,\Lambda)$.
\item
If every constructible subset of $X$ has locally finitely many irreducible components, then the natural t-structure on $\D(X,\Lambda)$ restricts to a t-structure on $\D_\cons(X,\Lambda)$.
\end{itemize}
\xtheo

The topological condition on $X$ is satisfied if $X$ is locally Noetherian. 
According to \thref{not.t.structure.rema}, \emph{some} topological assumption on $X$ is necessary for $\Dlis(X, \Ql)$ to be stable under truncation.
The following result allows for the comparison with the notions of lisse and constructible sheaves as in \cite{Deligne:Cohomologie, SGA5, Deligne:Weil2, Ekedahl:Adic, BhattScholze:ProEtale, StacksProject, GaitsgoryLurie:Weil}:

\theo
\thlabel{summary.sheaves}
Let $X$ be a scheme, and let $\Lambda$ be a condensed ring. 
\begin{enumerate}
  \item (\thref{discrete.comparison.prop})
  \label{item--summary.sheaves.1}
  For a discrete topological ring $\Lambda$, the \ii-category $\D_\lis(X, \Lambda)$ is equivalent to the full subcategory of $\D(X_\et,\Lambda)$ of complexes that are étale-locally perfect complexes of $\Lambda_*$-modules.
Consequently, $\D_\cons(X, \Lambda)$ is equivalent to the resulting category of étale constructible sheaves of $\Lambda$-modules.

  \item
  (\thref{limit.coefficients.lemm})
    \label{item--summary.sheaves.2}
Assume that $\Lambda = \lim \Lambda_n$ is a sequential limit of condensed rings such that the transition maps are surjective with locally nilpotent kernel.
Then there are natural equivalences 
\begin{align*}
  \D_\lis(X,\Lambda)\stackrel\cong \lr & \ \lim \D_\lis(X,\Lambda_n) \text{ and } \\
\D_\cons(X,\Lambda)\stackrel\cong \lr & \ \lim \D_\cons(X,\Lambda_n).
\end{align*}
 
  \item
  (\thref{localization.lisse.cons.lemm})
    \label{item--summary.sheaves.3}
 Assume that $T \subset \Lambda_*$ is a multiplicatively closed subset.
If $X$ is quasi-compact and quasi-separated (qcqs), then the natural functors
\begin{align*}
\D_\lis(X,\Lambda)\t_{\Perf_{\Lambda_*}}\Perf_{T^{-1}\Lambda_*} \r & \ \D_\lis\left(X,T^{-1}\Lambda\right) \text { and } \\
\D_\cons(X,\Lambda)\t_{\Perf_{\Lambda_*}}\Perf_{T^{-1}\Lambda_*} \r & \ \D_\cons\left(X,T^{-1}\Lambda\right)
\end{align*}
are fully faithful.

  \item
    \label{item--summary.sheaves.4}
  (\thref{colimit.coefficients.lemm})
Assume that  $\Lambda = \colim \Lambda_i$ is a filtered colimit of condensed rings. 
If $X$ is qcqs, there are natural equivalences
\begin{align*}
\colim \D_\lis(X,\Lambda_i)\stackrel \cong \lr & \ \D_\lis(X,\Lambda) \text{ and } \\
\colim \D_\cons(X,\Lambda_i)\stackrel \cong \lr & \ \D_\cons(X,\Lambda).
\end{align*}

  \item
    \label{item--summary.sheaves.5}
  (\thref{locally.constant.prop}) 
  Assume that $X$ has locally a finite number of irreducible components. 
  Then $M\in \D(X,\Lambda)$ is lisse if and only if $M$ is locally on $X_\proet$ isomorphic to $\underline{N}\t_{\underline{\Lambda_*}}\Lambda_X$ for some perfect complex of $\Lambda_*$-modules $N$, where the underline denotes the associated constant sheaf on $X_\proet$.
\end{enumerate}
\xtheo

Part \refit{summary.sheaves.2} applies to adic topological rings $\Lambda=\lim\Lambda/I^n$ such as the $\ell$-adic integers $\bbZ_\ell=\lim\bbZ/\ell^n$.
The key computation for the \ii-categories of perfect complexes is \cite[Lemma 4.2]{Bhatt:Algebraization}.
Together with \refit{summary.sheaves.1}, this also shows that $\D_\cons(X,\Lambda)$ is equivalent to the full subcategory of $\D(X,\Lambda)$ of $I$-adically complete objects $M$ such that $M\t \Lambda/I$ is étale constructible as in \cite[Section 6]{BhattScholze:ProEtale}.  
For algebraic field extensions $E\supset \bbQ_\ell$, one easily deduces from \refit{summary.sheaves.5} that the categories $\D_\cons(X,E)$, $\D_\cons(X,\calO_E)$ agree with the categories defined in \cite[Definition 6.8.8]{BhattScholze:ProEtale} whenever $X$ is topologically Noetherian. 
See, however, \thref{fun.example} for a lisse sheaf on a profinite set that is not pro\'etale-locally perfect-constant.
In conclusion, the above result extends the previously known approaches.

The functor in \refit{summary.sheaves.3} is not an equivalence for lisse sheaves in general.  
This relates to the difference between the \'etale versus the pro\'etale fundamental group, see \thref{etale.vs.proetale.fundamental.group} and also \thref{classical.comparison.lemm} for a positive result for constructible sheaves.
Part \refit{summary.sheaves.4} gives the comparison of the \ii-category of constructible sheaves with coefficients in $\bar\bbQ_\ell=\colim E$ or $\bar\bbZ_\ell=\colim \calO_E$ where the colimits run through finite field extensions $E\supset \bbQ_\ell$. 



Another application of the formalism presented in this work is the realization of  ind-constructible sheaf contexts as categories of sheaves. In algebraic geometry and geometric representation theory, it has been custumary and sometimes necessary to consider categories of inductive systems of constructible (or lisse) complexes \cite{Lafforgue:Chtoucas,Gaitsgory.Semi.infinite.I, AGKRRV.Local.Systems.Restricted.Langlands}.  For qcqs schemes $X$ of finite $\Lambda$-cohomological dimension (see \refsect{ind.cons.sheaves} for details), we introduce the following notion:
\defi 
\thlabel{ind.lis.cons.defi.intro}
A sheaf $M\in \D(X,\Lambda)$ is called \textit{ind-lisse}, respectively \textit{ind-constructible} if it is equivalent to a filtered colimit of lisse, respectively constructible sheaves.
\xdefi

We denote by $\D_\indlis (X,\Lambda)\subset \D_\indcons(X,\Lambda)$ the resulting full subcategories of $\D(X,\Lambda)$.
These \ii-categories satisfy \'etale descent, see \thref{descent.ind.cons}.
Examples of pairs $(X,\Lambda)$ satisfying the cohomological finiteness assumption include schemes $X$ of finite type over finite and separably closed fields with coefficients $\Lambda$ being a discrete torsion ring, an algebraic field extension $E\supset \bbQ_\ell$, or its ring of integers $\calO_E$, see \thref{Artin.vanishing}. For such pairs, we can realize categories of inductive systems as full subcategories of sheaves on the pro\'{e}tale site.

\prop[\thref{compact.objects}, \thref{Gaitsgory.Lurie.comparison}]
\thlabel{compact.objects.intro}
Let $X$ be a qcqs scheme of finite $\Lambda$-cohomological dimension. 
Then an object $M\in \D_\indcons(X,\Lambda)$ is compact if and only if $M$ is constructible, and likewise for (ind-)lisse sheaves. 
Consequently, the Ind-completion functor induces equivalences
\begin{align*}
\Ind\big(\D_\lis(X,\Lambda)\big)\stackrel\cong \lr & \ \D_\indlis(X,\Lambda) \text{ and } \\ 
\Ind\big(\D_\cons(X,\Lambda)\big)\stackrel\cong \lr & \ \D_\indcons(X,\Lambda).
\end{align*}
\xprop

Combined with \thref{summary.sheaves} \refit{summary.sheaves.1} and \cite[Section 6.4]{BhattScholze:ProEtale}, the proposition shows that $\D_\indcons(X,\Lambda)\cong \D(X_\et,\Lambda)$ for discrete rings $\Lambda$ and qcqs schemes $X$ of finite $\Lambda$-cohomological dimension.  

\rema
Another motivation for this work is the $\ell$-adic realization functor for motives. 
As explained in \cite{CisinskiDeglise:Etale}, this functor is in essence a completion functor. 
It seems an interesting question whether it can be expressed as a scalar extension functor for a yet-to-be-defined category of motives with condensed coefficients.
\xrema

\subsection*{Acknowledgements} 
We thank K\c{e}stutis \v{C}esnavi\v{c}ius, Ofer Gabber, Peter Scholze and Evgenij Vechtomov for helpful conversations and email exchanges.
The second named author (T.\ R.) thanks all participants of the GAUS seminar on Condensed Mathematics in the winter term 2020/2021. 
We heartily thank an anonymous referee for their extremely insightful report on our paper.

\section{Prelude on \ii-categories}
\label{sect--prelude}

Throughout this section, $\Lambda$ denotes a unital, commutative ring. 
We briefly collect some notation pertaining to \ii-categories from \cite{Lurie:HA,Lurie:Higher}.
As in \cite[Section 5.5.3]{Lurie:Higher}, $\PrL$ denotes the \ii-category of presentable \ii-categories with colimit-preserving functors. 
It contains the subcategory $\PrSt \subset \PrL$ consisting of stable \ii-categories.
Below, we also use the Ind-completion functor (defined on the \ii-category
of idempotent complete small stable \ii-categories and exact functors, 
taking values in compactly generated presentable \ii-categories) and the functor forgetting the compact generatedness \cite[Lemma~5.3.2.9]{Lurie:HA}
\eqn\Cat_\infty^{\perf} \stackrel [\cong]{\Ind} \lr \PrSt_{\omega} \lr \PrSt.
\label{eqn--CatExIdem.etc}
\xeqn

\subsection{Monoidal aspects}
\label{sect--monoidal.aspects}

The \ii-category $\PrL$ carries the Lurie tensor product \cite[Section 4.8.1]{Lurie:HA}.
This tensor product induces one on the full subcategory $\PrSt \subset \PrL$ consisting of stable \ii-categories \cite[Proposition~4.8.2.18]{Lurie:HA}.
For our commutative ring $\Lambda$, the \ii-category $\Mod_\Lambda$ of chain complexes of $\Lambda$-modules, up to quasi-isomorphism, is a commutative monoid in $\PrSt$ with respect to this tensor product.
This structure includes, in particular, the existence of a functor
$$\Mod_\Lambda \x \Mod_\Lambda \r \Mod_\Lambda$$
which, after passing to the homotopy categories is the classical \emph{derived} tensor product on the unbounded derived \ii-category of $\Lambda$-modules.

We define $\PrStL$ to be the category of modules, in $\PrSt$, over $\Mod_\Lambda$.
This \ii-category is denoted $\DGCat_{\cont, \Lambda}$ in \cite[Chapter 1, Section 10.3]{GaitsgoryRozenblyum:StudyI}, and we will freely use results from there.
Noting that modules over $\Mod_\Lambda$ are in particular modules over $\Sp$, the \ii-category of spectra, this can be described as the \ii-category consisting of \emph{stable} presentable \ii-categories together with a $\Lambda$-linear structure, and such that functors are continuous and $\Lambda$-linear.
Therefore $\PrStL$ carries a symmetric monoidal structure, whose unit is $\Mod_\Lambda$.

In order to express monoidal properties of \ii-categories consisting, say, of bounded complexes, 
recall that the tensor product of compactly generated \ii-categories are again compactly generated \cite[Lemma~5.3.2.11]{Lurie:HA}.
In addition, if $f^*, g^*$ are functors between such \ii-categories preserving compact objects, their right adjoints $f_*$ and $g_*$ preserve filtered colimits, and hence so does $f_* \t g_*$, which is the right adjoint of $f^* \t g^*$. Thus the latter preserves compact objects.
The symmetric monoidal structure on $\PrSt$ therefore restricts to one on $\PrSt_\omega$.
By \cite[Corollary~4.8.1.4]{Lurie:HA} or \cite[Proposition~4.4]{BenZviFrancisNadler:Integral} there is a symmetric monoidal structure on $\Cat_\infty^{\perf}$ characterized by
$$D_1 \t D_2 \defined \big(\Ind(D_1) \t \Ind(D_2)\big)^\omega,\eqlabel{Cat.idem.monoidal}$$
that is, the compact objects in the Lurie tensor product of the Ind-completions. 
With respect to these monoidal structures, both functors in \refeq{CatExIdem.etc} are symmetric monoidal.

The subcategory of compact objects in $\Mod_\Lambda$ is given by perfect complexes of $\Lambda$-modules \cite[Proposition~7.2.4.2.]{Lurie:HA}.
It is denoted $\Perf_\Lambda$.
Under the equivalence in \refeq{CatExIdem.etc}, the \ii-category $\Perf_\Lambda \in \Cat_\infty^\perf$ corresponds to $\Mod_\Lambda$. Moreover, $\Perf_\Lambda$ is a commutative monoid in $\Cat_\infty^{\perf}$, so that we can consider its category of modules, denoted as
$\Cat_{\infty, \Lambda}^{\perf}$.
This \ii-category inherits a symmetric monoidal structure denoted by $D_1 \t_{\Perf_\Lambda} D_2$.

Any stable \ii-category $D$ is canonically enriched over the \ii-category of spectra $\Sp$.
We write $\Hom_D(-, -)$ for the mapping spectrum.
Any object in $\PrStL$ is canonically enriched over $\Mod_\Lambda$, so that we refer to $\Hom_D(-, -) \in \Mod_\Lambda$ as the mapping complex.
For example, for $M, N\in \Mod_\Lambda$, then $\Hom_{\Mod_\Lambda}(M, N)$ is commonly also denoted by $\RHom(M, N)$.
Its $n$-th cohomology is the Hom-group $\Hom(M, N[n])$ in the classical derived category.

\subsection{Limits and filtered colimits}
\label{sect--limits.filtered.colimits}
Throughout, we freely use general facts about the forgetful functors 
$$\Cat_\infty^\perf \subset \CatEx_\infty\subset \Cat_\infty$$ 
between the \ii-categories of small stable idempotent complete, respectively small stable, respectively arbitrary small \ii-categories, together with exact, respectively exact, respectively all functors.
Recall that any functor between idempotent complete \ii-categories automatically preserves retracts or, equivalently, colimits indexed by $\Idem$.

\lemm
All three \ii-categories above have small limits and filtered colimits and both inclusions preserve these.
\xlemm

\pf 
See \cite[Theorem 1.1.4.4, Proposition 1.1.4.6]{Lurie:HA} for the claims concerning the latter functor.
For the former, we use that for any fully faithful right adjoint $D \subset C$ (such as $\iota: \Cat_\infty^\perf \subset \CatEx_\infty$), the inclusion creates all limits that $C$ admits.
Filtered colimits in $\Cat_\infty^\perf$ are computed by taking the idempotent completion of the colimit in $\CatEx_\infty$. According to \cite[Corollary~4.4.5.21]{Lurie:Higher}, however, that filtered colimit is already idempotent complete, showing that the inclusion $\iota$ also preserves filtered colimits.
\xpf

\lemm
\thlabel{perfect.complexes.colimit.lemm}
Let $\Lambda$ be a ring.
\begin{enumerate}
\item
\label{item--perfect.complexes.colimit.lemm.1}
If $\Lambda = \colim \Lambda_i$ is a filtered colimit of rings, then the natural functor
\[
\colim \Perf_{\Lambda_i} \lr \Perf_\Lambda
\]
is an equivalence. Here the transition functors are given by $(\str) \t_{\Lambda_i} \Lambda_j$ for $j \geq i$.


\item
\label{item--perfect.complexes.colimit.lemm.3}
Let $\Lambda=\lim_{i\geq 1} \Lambda_i$ be a sequential limit of rings such that all transition maps $\Lambda_{i+1}\r \Lambda_i$ are surjective with locally nilpotent kernel. 
Then the natural functor 
\[
\Perf_\Lambda \lr \lim\Perf_{\Lambda_i}, \;\; M\mapsto (M\t_\Lambda\Lambda_i)_{i\geq 1}
\]
is an equivalence. 
In addition, the functor $\Perf_\Lambda \r\Perf_{\Lambda_1}, M\mapsto M\t_\Lambda\Lambda_1$ is conservative. 
\end{enumerate}
\xlemm
\pf
In part (1), the full faithfulness follows from standard $\t$-$\Hom$-adjunctions using that perfect complexes are dualizable.
Since $\Lambda$ lies in the essential image, the functor is an equivalence given that both \ii-categories are in $\Cat_\infty^\perf$. 



Part (2) is \cite[Lemma 4.2]{Bhatt:Algebraization}: for full faithfulness, we note that the functor $M\t_{\Lambda}(\str)\cong \underline\Hom_{\Mod_\Lambda}(M^\vee,\str)$ commutes with limits so that $M\cong \lim M\t_{\Lambda}\Lambda_i$ for any $M\in \Perf_\Lambda$.
For essential surjectivity, we use that a quasi-inverse of the functor is provided by $\{M_i\}\mapsto \lim M_i$, see \StP{0CQG}. 

For the final statement, it remains to prove that the functor $\Perf_{\Lambda_i} \r\Perf_{\Lambda_1}, M\mapsto M\t_{\Lambda_i}\Lambda_1$ is conservative for any $i\geq 1$. 
So let $M\in \Perf_{\Lambda_i}$ such that $M\t_{\Lambda_i}\Lambda_1\simeq 0$.
We choose a bounded complex $(M^{a}\r \dots\r M^b)$ of finitely generated projective $\Lambda_i$-modules representing $M$.
We need to show that the complex is exact.
The derived tensor product $M\t_{\Lambda_i}\Lambda_1$ is represented by $(M^{a}\t_{\Lambda_i}\Lambda_1\r \dots\r M^b\t_{\Lambda_i}\Lambda_1)$.
As the complex is exact by assumption, the map $M^{b-1}\t_{\Lambda_i}\Lambda_1\r M^b\t_{\Lambda_i}\Lambda_1$ is surjective. 
Nakayama's lemma shows that $M^{b-1}\r M^b$ is surjective as well, using that $\ker(\Lambda_i\r \Lambda_1)$ is generated by nilpotent elements. 
By projectivity, this map splits so that $M^{b-1}\simeq M^b\oplus N^{b-1}$. 
We are reduced to show that the resulting complex $(M^a\r \dots\r M^{b-2} \r N^{ b-1})$ is exact. 
Continuing by induction on the length $b-a$ implies our claim. 
\xpf

Write $\Cat_n \subset \Cat_\infty$ for the full subcategory spanned by $n$-categories \cite[Section 2.3.4]{Lurie:Higher}, and similarly $\Cat_n(\Idem) \subset \Cat_\infty(\Idem)$ for the full subcategory consisting of idempotent complete $n$-categories.
Recall that $\Cat_n, \Cat_n(\Idem)$ admit limits and filtered colimits and that the above inclusions preserve these.
Let us record the following lemma from \cite[Lemma~3.7]{BhattMathew:Arc} for later use:

\lemm
\thlabel{filtered.colimits.limits.lemm}
In $\Cat_n$ and $\Cat_n(\Idem)$, filtered colimits commute with totalizations (=$\Delta$-indexed limits).
\xlemm
\pf
By \cite[Lemma 3.7, Example 3.6]{BhattMathew:Arc}, this holds true for $\Cat_n$.
Now use that $\Cat_n(\Idem) \subset \Cat_n$ preserves limits and filtered colimits. 

\xpf


\section{Definitions and basic facts}
\label{sect--generalities.proetale}

In this section, we start out setting up the basics of a theory of lisse (French for smooth) and constructible sheaves on schemes, with coefficients in a condensed (unital, commutative) ring $\Lambda$.

Throughout the paper, following \cite{BhattScholze:ProEtale}, we fix a strong limit cardinal $\kappa$. All schemes $X$ appearing in this paper are required to have $|X| < \kappa$. (Alternatively, one may avoid such a choice using the set-theoretic conventions in \cite[§1.2]{BarwickHaine:Pyknotic}.)

For a scheme $X$, we denote by $X_\proet$ its pro\'etale site introduced in \cite[\S4]{BhattScholze:ProEtale}. 
Thus, $X_\proet$ is the category of weakly étale schemes over $X$ with covers given by families $\{U_i\r U\}$ of maps in $X_\proet$ that are fpqc-coverings, that is, any open affine in $U$ is mapped onto by an open affine in $\sqcup_i U_i$.
Being a major difference with étale site $X_\et$, any scheme $X$ admits a pro\'etale cover by w-contractible affine schemes, that is, by affine schemes $U$ such that any weakly étale surjection $V\r U$ splits.

Also, for any profinite set $S$, we denote by $S_\proet$ the category of profinite sets over $S$ with covers given by finite families of jointly surjective maps.
The 
w-contractible objects in $S_\proet$ are the extremally disconnected profinite sets. 
We use this in the special case where $S=*$ is the singleton: For any scheme $X$, there is the map of sites 
\[
p_X\co X_\proet\lr *_\proet
\] 
given by the limit-preserving functor (in the opposite direction) $S=\lim S_i \mapsto \lim X^{S_i}=:X\x S$.

Throughout this section, $\Lambda$ is a condensed (unital, commutative) ring $\Lambda$, that is, a sheaf of rings on $*_\proet$, see \cite{ClausenScholze:CondensedMathematics}.
We write $\Lambda_X:=p_X^{-1}\Lambda$ for the corresponding sheaf of rings on $X_\proet$.
The condensed rings considered in this paper arise as follows:

\exam
\thlabel{topological.rings.exam}
Every T1-topological ring $\Lambda$ induces sheaf of rings, also denoted by $\Lambda$, on $*_\proet$:
\[
 S\longmapsto \Maps_\cont(S,\Lambda),
\]
that is, continuous maps from the profinite set $S$ to $\Lambda$.
In fact, we have a functor from the category of T1-topological rings into the category of condensed rings.
This functor is fully faithful when restricted to those $\Lambda$ whose underlying topological space is compactly generated.

In fact, for any scheme or profinite set $X$, the sheaf $\Lambda_X$ is nothing but the sheaf of rings on $X_\proet$ given by 
\[
U\longmapsto \Maps_\cont(U,\Lambda), 
\]
that is, continuous maps from the underlying topological space of $U$ into $\Lambda$: any continuous map $U\to \Lambda$ from a w-contractible affine $U\in X_\proet$ factors uniquely through $U\r \pi_0U$ because $\Lambda$ is T1. 
In \cite[Lemma 4.2.12]{BhattScholze:ProEtale}, this sheaf is denoted $\calF_\Lambda$.

Favorable examples of T1-topological rings $\Lambda$ include all discrete rings, adic rings,
algebraic field extensions $E\supset \bbQ_\ell$ equipped with the colimit topology, their open subrings of integers $\calO_E$, but also the ring of adeles $\bbA_{K}^T$ prime to some finite set of places $T$ in some number field $K$.
\xexam

For any scheme $X$, we denote by $\D(X,\Lambda)$ the derived \ii-category of the abelian category of sheaves of $\Lambda_X$-modules on $X_\proet$. 
This is a presentable stable \ii-category which is symmetric monoidal and closed.
The monoidal structure is denoted by $\str\otimes_{\Lambda_X}\str$ and the inner homomorphisms by $\underline \Hom_{\Lambda_X}(\str,\str)$, so that 
\[
\Hom_{\Lambda_X}(\str,\str)=\Hom_{\D(X,\Lambda)}(\str,\str)=\RG\big(X,\underline \Hom_{\Lambda_X}(\str,\str)\big)
\]
is the mapping complex. 

For any $n\in \bbZ$, the truncations $\tau^{\geq n}$, $\tau^{\leq n}$ in the standard t-structure on $\D(X,\Lambda)$ induce (in cohomological notation) adjunctions
\[
\tau^{\geq n} : \D(X,\Lambda) \leftrightarrows \D^{\geq n}(X,\Lambda) : \on{incl}, \;\;\;\;  \on{incl}: \D^{\leq n}(X,\Lambda) \leftrightarrows \D(X,\Lambda) : \tau^{\leq n},
\]
where $\on{incl}$ is the inclusion of the respective subcategories. 
The subcategory $\D^{\le 0}(X, \Lambda)$ is preserved under $\str \t_{\Lambda_X} \str$.
As $X_\proet$ is locally weakly contractible \cite[Proposition 4.2.8]{BhattScholze:ProEtale}, the t-structure on $\D(X,\Lambda)$ is left-complete (equivalently, Postnikov towers converge in the associated hypercompleted \ii-topos) and the \ii-category $\D(X,\Lambda)$ is compactly generated, see \cite[Proposition 3.2.3]{BhattScholze:ProEtale}.
A family of compact generators is given by the objects $\Lambda_X[U]\in \D(X,\Lambda)$, for $U\in X_\proet$ w-contractible affine, corepresenting the functor $\RG(U,\str)$. 
We remark that $\D(X,\Lambda)$ is equivalent to the \ii-category of $\Lambda_X$-modules on the hypercompleted \ii-topos associated with $X_\proet$ by \cite[Theorem 2.1.2.2, Definition 2.1.0.1]{Lurie:SAG}. 
The above also applies to any profinite set $S$ and the \ii-category $\D(S,\Lambda)$ in place of $X$ and $\D(X,\Lambda)$ (in fact, this is just a special case).

If $f\co Y\r X$ is any morphism of schemes, then $\Lambda_Y=f^{-1}\Lambda_X$.
It is formal to check that the ordinary pullback, respectively pushforward of sheaves induces an adjunction
\[
f^*=f^{-1} : \D(X,\Lambda)\leftrightarrows \D(Y,\Lambda) : f_*,
\]
where $f^*$ is exact, t-exact and symmetric monoidal. 
Similarly, if $\Lambda\r \Lambda'$ is a morphism of condensed rings, then the forgetful functor $\D(X,\Lambda')\r \D(X,\Lambda)$ admits a symmetric monoidal left adjoint 
\[
\D(X,\Lambda)\r \D(X,\Lambda'),\;\; M\mapsto M\otimes_{\Lambda_X}\Lambda_X'.
\]

Let us denote by $\Gamma(X,\Lambda)$ the (underived) global sections of $\Lambda_X$ viewed as a ring. 
Then there is the functor
\eqn
\label{constant.objects.functor}
\Mod_{\Gamma(X,\Lambda)} \r \D(X,\Lambda), \;\; M\mapsto M_X,
\xeqn
characterized in $\PrSt$ as the colimit-preserving extension of $\Gamma(X,\Lambda)\mapsto \Lambda_X$.
Explicitly, $M_X$ is the hypersheaf associated with the presheaf on $X_\proet$ given by $U\mapsto M\otimes_{\Gamma(X,\Lambda)}\RG(U,\Lambda)$.
This functor is symmetric monoidal and makes $\D(X,\Lambda)$ a commutative algebra object in $\PrSt_{\Gamma(X,\Lambda)}$.
The sections of $M_X$ on w-contractible qcqs $U\in X_\proet$ are computed as 
\eqn
\label{global.sections.dualizable}
\RG(U,M_X)\cong M\otimes_{\Gamma(X,\Lambda)}\Gamma(U,\Lambda).
\xeqn
In particular, if $f\co Y\r X$ lies in $X_\proet$, then $f^*M_X\cong \big(M\t_{\Gamma(X,\Lambda)}\Gamma(Y,\Lambda)\big)_Y$.

\rema
Here is an equivalent way of defining the functor $M\mapsto M_X$ in \eqref{constant.objects.functor}:
For any scheme or profinite set $X$, there is a natural map $\underline {\Gamma(X,\Lambda)}\r \Lambda_X$ of sheaf of rings on $X_\proet$ where $\underline {\Gamma(X,\Lambda)}$ denotes the constant sheaf associated with the discrete ring $\Gamma(X,\Lambda)$.
Then the functor \eqref{constant.objects.functor} is equivalent to the functor 
\[
M\longmapsto \underline M \t_{\underline {\Gamma(X,\Lambda)}}\Lambda_X,
\]
where $\underline M$ denotes the constant sheaf.
As an example, let $X=*$, $\Lambda=\bbQ_\ell$ (see \thref{topological.rings.exam}) and $M\in \Mod_{\bbQ_\ell}^\heartsuit$ a $\bbQ_\ell$-vector space.
Then, loosely speaking, the above functor equips $M$ with the relatively discrete topology. 
More precisely, writing $M=\colim \bbQ_\ell^I$ as an increasing union of finite-dimensional vector spaces, we take the product topology on $\bbQ_\ell^I$ and the colimit topology on $M$.
\xrema

Recall that a subset of a qcqs topological space is called \emph{constructible} if it is a finite Boolean combination of quasi-compact open subsets. 
Also, recall the notion of dualizable objects in symmetric monoidal categories \cite[Definition 4.6.1.1, Remark 4.6.1.12]{Lurie:HA}.

\defi
\thlabel{lisse.constructible.defi}
Let $X$ be a scheme or a profinite set, and $\Lambda$ a condensed ring. 
\begin{enumerate}
\item 
A sheaf $M\in \D(X,\Lambda)$ is called {\it lisse} if it is dualizable.

\item 
A sheaf $M\in \D(X,\Lambda)$ is called {\it constructible} if, for any open affine $U\subset X$, there exists a finite subdivision of $U$ into constructible locally closed subschemes $U_i\subset U$ such that $M|_{U_i}$ is lisse.
\end{enumerate}
\xdefi

If $X$ is qcqs (=quasi-separated and quasi-compact) and $M\in \D(X,\Lambda)$ constructible, then there is a finite subdivision of $X$ into constructible locally closed subschemes $X_i\subset X$ such that $M|_{X_i}$ is lisse.
The argument is purely topological and the same as in \StP{095E}.

The full subcategories of $\D(X,\Lambda)$ of lisse, respectively constructible $\Lambda$-sheaves are denoted by
\[
\D_\lis(X,\Lambda) \;\subset\; \D_\cons(X,\Lambda).
\]
Both \ii-categories are naturally commutative algebra objects in $\Cat_{\infty, \Gamma(X,\Lambda)}^\perf$, that is, idempotent complete stable $\Gamma(X,\Lambda)$-linear symmetric monoidal \ii-categories.

If $f\co Y\r X$ is any map of schemes, then the pullback $f^*\co \D(X,\Lambda)\r \D(Y,\Lambda)$ preserves lisse, respectively constructible sheaves and hence induces functors
\[
f^*\co \D_\lis(X,\Lambda)\r \D_\lis(Y,\Lambda),\;\;\;\; f^*\co \D_\cons(X,\Lambda)\r \D_\cons(Y,\Lambda).
\]
For lisse sheaves, this follows from the monoidality of $f^*$.
For constructible sheaves, one additionally reduces to the case of affine schemes so that $f$ induces a spectral map on the underlying topological spaces, and thus is continuous in the constructible topology, see \StP{0A2S}.

If $\Lambda\r \Lambda'$ is a map of condensed rings, then the functor $(\str)\t_{\Lambda_X}\Lambda_X'$ preserves lisse, respectively constructible sheaves and hence induces functors
\[
\D_\lis(X,\Lambda)\r \D_\lis(X,\Lambda'),\;\;\;\; \D_\cons(X,\Lambda)\r \D_\cons(X,\Lambda').
\]

For any constructible closed immersion $i\co Z\hr X$ with open complement $j\co U\hr X$, we have adjunctions
\eqn
\label{adjoints.closed.open.immersion.eq}
j_! : \D(U,\Lambda)\leftrightarrows \D(X,\Lambda) : j^*,\;\;\;  i_*: \D(Z,\Lambda)\leftrightarrows \D(X,\Lambda) : i^!,
\xeqn
fitting in fiber sequences $j_!j^*\r \id \r i_*i^*$ and $i_*i^!\r \id \r j_*j^*$, see \cite[\S6.1]{BhattScholze:ProEtale}.
The functors $j_*, j_!,i_*$ are fully faithful and satisfy the usual formulas
\eqn
\label{formulas.closed.open.immersion.eq}
i^*i_* \simeq j^*j_*\simeq j^*j_!\simeq \id, \;\;\;\; j^*i_* \simeq i^*j_!\simeq 0.
\xeqn

\lemm
\thlabel{t.exactness.lemm}
In the above situation, the functors $i_*, j_!$ are t-exact and preserve the full subcategories of constructible sheaves.
\xlemm
\pf
By \cite[Lemma 6.2.1 (1) + (2)]{BhattScholze:ProEtale}, the functors $i_*,j_!$ induce equivalences onto the full subcategory of $\D(X_\proet,\Lambda)$ spanned by objects supported on $Z$, respectively $U$.
Their inverses are given by $i^*,j^*$ which are clearly $t$-exact, 
hence so are the functors $i_*,j_!$.
Using the formulas \eqref{formulas.closed.open.immersion.eq} it is clear that $i_*,j_!$ preserve constructibility. 
\xpf

We use the following terminology throughout:

\defi
\thlabel{perfect.lc}
 A sheaf $M \in \D(X, \Lambda)$ is called \emph{(perfect-)constant} if $M\simeq\underline{N}\t_{\underline{\Lambda_*}}\Lambda_X$ for some (perfect) complex of $\Lambda_*$-modules $N$, where $\Lambda_*=\Gamma(*,\Lambda)$ is the underlying ring.
 It is called {\it (pro-)\'etale-locally (perfect-)constant} if it is so locally on $X_{\text{(pro-)\'et}}$.
\xdefi

Any pro\'etale-locally perfect-constant sheaf is lisse. 
The converse holds for discrete coefficient rings $\Lambda$ (\thref{locally.constant.lemm}) in which case $\underline{\Lambda_*}\cong\Lambda_X$. 
It also holds for schemes having locally finitely many irreducible components, see \thref{locally.constant.prop}.
The following example of a lisse sheaf, which is based on \cite[Example 6.6.12]{BhattScholze:ProEtale}, shows however that lisse sheaves on profinite sets do not have such a simple description:

\exam
\thlabel{fun.example}
Let $S=\hat\bbZ=\lim_m \bbZ/m$ be the profinite completion of the integers viewed as a profinite set. 
Take $\Lambda=\hat\bbZ$ viewed as a profinite ring. 
Then the endomorphisms of $\Lambda_S$ in $\D(S,\Lambda)^\heartsuit$ are computed as
\[
\H^0(S,\Lambda)=\Maps_\cont(S,\hat \bbZ)=\Maps_\cont(\hat \bbZ,\hat \bbZ).
\]
The constant map $f\equiv s \in \hat\bbZ$ corresponds to the endomorphism of $\Lambda_S$ given by multiplication with the scalar $s$.
By contrast, if $f\co \Lambda_S\r \Lambda_S$ corresponds to the identity in $\Maps_\cont(\hat \bbZ,\hat \bbZ)$, then its stalk $f_s$ at $s\in S$ is multiplication with $s$ viewed as element in $\hat\bbZ$.
The complex $(\Lambda_S \stackrel f \r \Lambda_S)\in \D_\lis(S,\Lambda)$ is \'etale-locally perfect-constant after each reduction modulo $m\neq 0$. 
However, the complex is not pro\'etale-locally perfect-constant as any cover $\{S_i\r S\}$ has a member $S_i\r S$ whose image is infinite. 
\xexam

\section{Descent properties}

Lisse and constructible sheaves interact nicely with proétale descent.
We first study lisse sheaves on w-contractible schemes, which may be thought of as a basis for the proétale topology of objects having cohomological dimension zero.
We then establish hyperdescent for $\Dlis$ and $\Dcons$ and draw some consequences.

\subsection{Lisse sheaves on w-contractible schemes}
Recall the functor $\Mod_{\Gamma(X,\Lambda)} \r \D(X,\Lambda), M\mapsto M_X$ from \eqref{constant.objects.functor}.

\lemm
\thlabel{dualizable.objects.contractible.lem}
Let $X$ be a w-contractible qcqs scheme, or an extremally disconnected profinite set.
\begin{enumerate}
\item 
\label{item--dualizable.objects.contractible.lem.1}
The functor \eqref{constant.objects.functor} induces an adjunction
$$(\str)_X : \Mod_{\Gamma(X, \Lambda)} \rightleftarrows \D(X, \Lambda) : \RG(X,\str).$$
Both adjoints are colimit-preserving and symmetric monoidal. 
In addition, $(\str)_X$ is fully faithful and $\RG(X,\str)$ is t-exact. 
\item
\label{item--dualizable.objects.contractible.lem.3}
The adjunction induces an equivalence on dualizable objects:
$$(\str)_X : \Perf_{\Gamma(X, \Lambda)}\cong \D_\lis(X, \Lambda) : \RG(X,\str)$$
\end{enumerate}
\xlemm

\pf
For \refit{dualizable.objects.contractible.lem.1}, we note that the spectra-valued functor $\RG(X,\str)\co \D(X, \Lambda)\r \Sp$ is right adjoint to the unique colimit-preserving functor $\Sp \r \D(X, \Lambda)$ mapping the sphere spectrum to $\Lambda_X$. 
Since $\Mod_{\RG(X,\Lambda)}$ is the category of modules over the monad associated to the adjunction, this induces an adjunction 
\[
\Mod_{\RG(X, \Lambda)} \rightleftarrows \D(X, \Lambda).
\]
Since $X$ is w-contractible qcqs, respectively extremally disconneced, it is w-contractible coherent in the topos-theoretic sense.
Thus $\RG(X,\str)$ is limit-preserving, t-exact and preserves colimits of uniformly bounded below diagrams, hence all colimits using the left-completeness of $\D(X, \Lambda)$.
In particular, $\RG(X,\Lambda)=\Gamma(X,\Lambda)$.
The unit of the adjunction $\id \r\RG(X,\str)\circ (\str)_X$ is an equivalence by \eqref{global.sections.dualizable} so that $(\str)_X$ is fully faithful.
Clearly, $(\str)_X$ is also colimit-preserving and symmetric monoidal.
For $M, N\in \D(X,\Lambda)$, their tensor product $M\t_{\Lambda_X}N$ is the \ii-sheafification of the presheaf $U\mapsto \RG(U,M)\t_{\RG(U,\Lambda)}\RG(U,N)$. 
Again, since $X$ is w-contractible coherent, the sheafification is unnecessary, so its global sections are equivalent to $\RG(X,M)\t_{\Gamma(X,\Lambda)}\RG(X,N)$. 

For \refit{dualizable.objects.contractible.lem.3}, we note that the adjunction restricts to an adjunction on dualizable objects by monoidality of both functors.
The counit of this adjunction $(\str)_X\circ \RG(X,\str)\r \id$ is an equivalence if and only if the functor
\[
\RG(X,\str)\co \D_\lis(X, \Lambda)\r \Perf_{\Gamma(X,\Lambda)}
\]
is fully faithful, see \cite[Proposition 5.2.7.4]{Lurie:Higher}.
For $M,N\in \D_\lis(X,\Lambda)$, this follows from $\underline \Hom_{\Lambda_X}(M,N)\cong N\t_{\Lambda_X} M^\vee$ upon applying the symmetric monoidal functor $\RG(X,\str)$.
\xpf

The following example, communicated to us by Peter Scholze, shows that the functor $(\str)_X$ is not t-exact in general. 

\exam 
\thlabel{not.t.exact.example}
Let $X=\beta\bbN$ be the Stone-\v{C}ech compactification of the natural numbers viewed as an extremally disconnected profinite set.
Let $\Lambda=\bbQ_\ell$ viewed as a condensed ring, see \thref{topological.rings.exam}.
The map $\bbN\r \bbQ_\ell$, $n\mapsto \ell^n$ uniquely extends to a continuous map $f\co \beta\bbN\r \bbQ_\ell$ by the universal property of $\beta$, that is, $f\in\Maps_\cont(\beta\bbN,\bbQ_\ell)=\Gamma(X,\Lambda)$. 
One checks that the complex $0\r \Gamma(X,\Lambda)\stackrel f \r \Gamma(X,\Lambda)$ is exact.
However, the induced complex on the level of sheaves is not exact because $f|_{\partial X}=0$, where $\partial X= \beta\bbN\backslash \bbN$ denotes the boundary.
\xexam

If $X$ is a qcqs scheme, then its underlying topological space is spectral \StP{094L}. 
Thus, the set of connected components $\pi_0X$ endowed with the quotient topology is a profinite space.
Any map $S\r \pi_0X$ of profinite sets can be written as profinite $\pi_0X$-sets $S=\lim S_i$ such that each $S_i\r \pi_0X$ is the base change of a map of finite sets, see \cite[Proof of Lemma 2.2.8]{BhattScholze:ProEtale}.
If we equip the topological space $X\x_{\pi_0X}S_i\r X$ with the sheaf of rings given by the pullback of the structure sheaf on $X$, then it is representable by an object of the Zariski site $X_\Zar$. 
The induced transition maps $X\x_{\pi_0X}S_j\r X\x_{\pi_0X}S_i$, $j\geq i$ are affine so that the limit
\[
X\x_{\pi_0X}S\defined \lim X\x_{\pi_0X}S_i
\]
exists in the category of $X$-schemes.
The functor $S\mapsto X\x_{\pi_0X}S$ from profinite $\pi_0X$-sets to $X$-schemes is limit-preserving and induces a map of sites
\begin{equation}
\label{connected.component.map.eq}
\pi_X\co X_\proet\lr (\pi_0X)_\proet,
\end{equation}
factorizing $p_X\co X_\proet\r *_\proet$.

\prop
\thlabel{connected.component.pullback}
Let $X$ be a w-contractible affine scheme.
\begin{enumerate}
\item
\label{item--connected.component.pullback.1}
The functor
\[
\pi_X^*\co \D( \pi_0X,\Lambda) \lr  \D( X,\Lambda)
\]
is fully faithful and commutes with the formation of inner homomorphisms. 

\item 
\label{item--connected.component.pullback.2}
A $\Lambda$-sheaf $M\in \D( X,\Lambda)$ lies in the essential image of $\pi_X^*$ if and only if for all maps $U\r V$ in $X_\proet$ between w-contractible affine schemes inducing isomorphisms $\pi_0U\cong\pi_0V$, the map
\[
\RG(V,M)\stackrel\cong \lr \RG(U,M)
\]
is an equivalence.

\item 
\label{item--connected.component.pullback.3}
The functor $\pi_X^*$ induces an equivalence 
$$\D_\lis(\pi_0 X, \Lambda) \stackrel \cong \lr \D_\lis(X, \Lambda).$$
\end{enumerate}
\xprop
\pf

We adjust the argument given in \cite[Lemma 4.2.13]{BhattScholze:ProEtale}  for the abelian categories. 
Abbreviate $\pi=\pi_X$, $\D(X)=\D(X,\Lambda)$ and $\D(\pi_0X)=\D(\pi_0X,\Lambda)$. 

For \refit{connected.component.pullback.1}, we show that the natural map $\id\r \pi_*\pi^*$ is an equivalence which formally implies the full faithfulness.
Any continuous $\pi_0X$-map $U\r S$ with affine $U\in X_\proet$, $S\in (\pi_0X)_\proet$ factors uniquely through $U\r \pi_0U$, since any profinite set is totally disconnected.  
Hence, if $M\in \D(\pi_0X)$, then $\pi^*M$ is the sheafification of the presheaf $U \mapsto \RG(\pi_0U,M)$.
In particular, if $U$ is also w-contractible, then we have an equivalence
\eqn
\label{evaluation.formula.eq}
\RG(U,\pi^*M) \cong \RG(\pi_0U,M).
\xeqn
In this case, $\pi_0U$ is extremally disconnected by \cite[Lemma 2.4.8]{BhattScholze:ProEtale}.
We apply these observations to show that the map $M\r \pi_*\pi^*M$ is an equivalence as follows. 
By evaluating at any extremally disconnected $S\in(\pi_0X)_\proet$ it suffices to show that the map
\[
\RG(S,M)\lr \RG\left(X\x_{\pi_0X}S,\pi^*M\right)
\]
is an equivalence. 
As $X\x_{\pi_0X}S\r X$ is an pro-(Zariski open) pro-finite map with $\pi_0(X\x_{\pi_0X}S)\cong S$ (by construction) we see that $X\x_{\pi_0X}S$ is w-contractible affine as well:
affine is clear; $X\r \pi_0X$ has a section $s\co \pi_0X\r X$ given by the closed points in $X$ by w-locality \cite[Lemma~2.1.4]{BhattScholze:ProEtale} so that 
\[
s\x_{\pi_0X} \id \co S=\pi_0X\x_{\pi_0X}S\r X\x_{\pi_0X}S 
\]
identifies $S$ with the closed points in $X\x_{\pi_0X}S$; finally, $X\x_{\pi_0X}S\r X$ induces an isomorphism on local rings which are therefore strictly Henselian at all closed points. 
This shows that $X\x_{\pi_0X}S$ is w-strictly local and hence w-contractible by \cite[Lemma 2.4.8]{BhattScholze:ProEtale} using that its set of connected components is $S$ (which is extremally disconnected).
This implies $M\cong \pi_*\pi^*M$.
The preservation of inner homomorphisms is immediate from the full faithfulness using \eqref{evaluation.formula.eq}.

For \refit{connected.component.pullback.2}, if $M\in \D(X)$ is equivalent to the $\pi^*$-pullback of some object in $D(\pi_0X)$, then it satisfies the desired condition by \eqref{evaluation.formula.eq}.
Conversely, assume that $M$ is localizing for maps $U\r V$ in $X_\proet$ of w-contractible affine schemes inducing an isomorphism on $\pi_0$.
We claim that the map $\pi^*\pi_*M\r M$ is an equivalence. 
Indeed, evaluating at some w-contractible affine $U\in X_\proet$ gives the map
\eqn
\label{evaluation.formula.eq2}
\RG(X\x_{\pi_0X}\pi_0U,M)\lr \RG(U,M)
\xeqn
induced from the canonical map $U\r X\x_{\pi_0X}\pi_0U$ over $\pi_0U$. 
One argues as in (1) above to see that $X\x_{\pi_0X}\pi_0U$ is w-contractible affine with space of components $\pi_0U$.  
Thus, \eqref{evaluation.formula.eq2} is an isomorphism by our assumption on $M$.
We conclude $M\cong \pi^*\pi_*M$.

For \refit{connected.component.pullback.3}, we note $\Lambda_X\cong \pi^*\Lambda_{\pi_0X}$ so that $\Gamma(X,\Lambda)=\Gamma(\pi_0X,\Lambda)$ by (1).
By \eqref{evaluation.formula.eq} the diagram
\[
\xymatrix{
 \D_\lis(\pi_0X) \ar[r]^{\pi^*} & \D_\lis(X) 
\\
\Perf_{\Gamma(\pi_0X,\Lambda)} \ar@{=}[r] \ar[u]^{\cong}_{(\str)_{\pi_0X}} & \Perf_{\Gamma(X,\Lambda)} \ar[u]^{\cong}_{(\str)_{X}}.
}
\]
commutes up to equivalence.
The vertical functors are equivalences by \thref{dualizable.objects.contractible.lem} \refit{dualizable.objects.contractible.lem.3}.
\xpf

The following corollary shows that lisse sheaves extend pro\'etale locally to small neighborhoods:

\coro
\thlabel{neighborhood.extension.coro}
Let $A$ be a ring Henselian along an ideal $I$. 
Let $i\co Z:=\Spec A/I\hr \Spec A=:X$ be the closed immersion induced by the quotient map $A\r A/I$.
\begin{enumerate}
\item 
\label{item--neighborhood.extension.coro.1}
The map $i$ induces an isomorphism $\pi_0Z\cong \pi_0X$. 

\item
\label{item--neighborhood.extension.coro.2}
The affine scheme $X$ is w-contractible if and only if $Z$ is w-contractible.
In this case, there is an equivalence
\[
i^*\co \D_\lis(X,\Lambda)\stackrel\cong\lr \D_\lis(Z,\Lambda).
\]
\end{enumerate}
\xcoro
\pf
For \refit{neighborhood.extension.coro.1}, the following argument was explained to us by K\c{e}stutis \v{C}esnavi\v{c}ius: 
Being profinite sets the map $\pi_0Z\r \pi_0X$ is obtained as a limit over (finer and finer) finite subdivisions of $Z$ and $X$ into clopen (=closed and open) subsets. 
By the unique lifting of idempotents along the quotient map $A\r A/I$ \StP{09XI}, these finite subdivisions match, so do the limits.  

For \refit{neighborhood.extension.coro.2}, we note that by \cite[Theorem 1.8, Lemma 2.2.9]{BhattScholze:ProEtale} an affine scheme is w-contractible if and only if it is w-strictly local and its profinite set of connected components is extremally disconnected.
Since $A$ is Henselian along $I$, the ring $A$ is w-strictly local if and only if $A/I$ is w-strictly local, see \cite[Lemma 2.2.13]{BhattScholze:ProEtale}. 
So the first statement in \refit{neighborhood.extension.coro.2} follows from \refit{neighborhood.extension.coro.1}. 
The second then follows, again using \refit{neighborhood.extension.coro.1}, from \thref{connected.component.pullback} \refit{connected.component.pullback.3}.
\xpf 

\subsection{Hyperdescent}
\label{sect--hyperdescent}
In this subsection, let $X$ be a scheme and $\Lambda$ a condensed ring.

\lemm
\thlabel{lisse.constructible.local.property.lemm}
The property of $\Lambda$-sheaves of being lisse, respectively constructible is local on $X_\proet$. 
\xlemm
\pf
It is enough to prove the following: if $X$ is affine and $j\co U\r X$ a w-contractible affine cover, then $M\in \D(X,\Lambda)$ is lisse, respectively constructible if and only if $j^*M$ is so. 
Since $j^*\co \D(X,\Lambda)\r \D(U,\Lambda)$ is monoidal, conservative and commutes with inner homomorphisms, the statement for the property ``lisse'' follows. 

Now assume that $j^*M$ is constructible. 
Since the $1$-topos of $X_\proet$ is generated by pro-étale affine objects \cite[Lemma 4.2.4]{BhattScholze:ProEtale}, we can assume that $U=\lim_iU_i\r X$ is a cofiltered limit of affine schemes $U_i\in X_\et$. 
If the stratification on $U$ witnessing the constructibility of $j^*M$ arises by pullback from $X$, then we are done using the case of lisse sheaves above. 
Following \cite[Lemmas 6.3.10, 6.3.13]{BhattScholze:ProEtale} we reduce to this situation in several steps.

First, each constructible subset of $U$ arises by pullback from some $U_i$. 
So the stratification witnessing the constructibility of $j^*M$ arises by pullback from some $U_i$.
We reduce to the case where $U=U_i\r X$ is an étale cover.

Next, stratifying $X$ by constructible locally closed subschemes $X_i\subset X$ such that the base change $U\x_XX_i\r X_i$ is finite étale \StP{03S0} we may assume that $U\r X$ is finite étale (after replacing $X$ by some $X_i$, and possibly a Zariski localization to preserve the affineness of $X$). 

Now writing $X=\lim_i X_i$ as a cofiltered inverse limit of finite type $\bbZ$-schemes the map $U\r X$ arises as the base change of some finite étale map $U_i\r X_i$.
The connected components of $X_i$ are open and closed.
After possibly replacing $X$ by a finite clopen cover we may assume that $X_i$ is connected.
Likewise, we may assume that $U_i$ is connected. 
Then we may replace $U_i\r X_i$ by its Galois closure and assume that $U_i\r X_i$ is a finite Galois cover with group $G=\Aut(U_i/X_i)$. 
Hence, we reduced to the case where $j\co U\r X$ is a $G$-torsor under some finite constant $X$-group scheme $G$. 

Finally, using the $G$-action on $U$ one easily constructs a finite subdivision of $U$ into $G$-equivariant constructible locally closed $U_i\subset U$ such that $j^*M|_{U_i}$ is lisse.
Clearly, these strata arise by pullback along the $G$-torsor $U\r X$.
This implies the constructibility of $M$ (again using the case of lisse sheaves above).
\xpf

\coro
\thlabel{preservation.constructibility}
If $j\co U\r X$ is quasi-compact \'etale (respectively, finite \'etale), then $j_!\co \D(U,\Lambda)\r \D(X,\Lambda)$ preserves the subcategories of constructible sheaves (respectively, lisse sheaves).
\xcoro
\pf
As in the proof of \thref{lisse.constructible.local.property.lemm}, one reduces to the finite \'etale case and further to the case of a $G$-torsor $j\co U\r X$ under some finite constant $X$-group scheme $G$.
Then $U\x_XU\cong G\x X$ which implies $j^*j_!M\cong\oplus_{g\in G} M$ for any $\D(X,\Lambda)$.
The corollary follows.  
\xpf

\coro 
\thlabel{lisse.constructible.hyperdescent.coro}
The functors $U\mapsto \D_\cons(U, \Lambda), \ \D_\lis(U, \Lambda)$ are hypersheaves of \ii-categories on $X_\proet$.
\xcoro
\pf 
We only spell out the constructible case, the one for lisse sheaves is identical. 
We first check that $U\mapsto \D_\cons (U,\Lambda)$ is a sheaf on $X_{\proet}$. Given an object $U\in X_\proet$ and an \'{e}tale cover $\{U_i \rightarrow U\}$ we denote by $\mathcal{U}$ the covering sieve generated by the maps $\{U_i \rightarrow U\}$. By \cite[Remark 2.1.0.5]{Lurie:SAG}, we know that the functor $U \mapsto \D(U,\Lambda)$ is a hypersheaf on $X_\proet$. 
Thus, we have an equivalence
\[
\D(U,\Lambda) \;\stackrel \cong\lr\; \lim_{V\in \mathcal{U}} \D(V,\Lambda). 
\]
As constructibility is preserved by pullback,
we have inclusions of full subcategories
\[
\D_\cons(U,\Lambda) \;\subset\; \lim_{V\in \mathcal{U}} \D_\cons (V,\Lambda) \;\subset\; \D(U,\Lambda).
\]
The essential image of the limit consists of objects $M\in \D(U,\Lambda)$ such that $M|_V$ is constructible for every $V\in \mathcal{U}$. 
In particular, $M|_{U_i}$ is constructible for every $U_i$ in the cover $\{U_i \rightarrow U\}$. Hence, $M\in \D_\cons(U,\Lambda)$ by \thref{lisse.constructible.local.property.lemm}.

Given the sheaf property, being a hypersheaf can be checked locally, we reduce to the case where $X$ is affine. 
In this case, the Grothendieck topology $X_\proet$ is finitary in the sense of \cite[Section A.3.2]{Lurie:SAG}. 
So by \cite[Proposition A.5.7.2]{Lurie:SAG} it is enough to show that for every hypercover $U_\bullet \r U$ with $U\in X_\proet$ qcqs, the natural functor 
    \[
    \D_\cons(U,\Lambda) \lr \tot\big(\D_\cons(U_\bullet,\Lambda)\big) := \lim_{[n]\in \Delta} \big(\D_\cons(U_{n},\Lambda)\big).
    \]
 is an equivalence.
Since $U\mapsto \D(U,\Lambda)$ is a hypersheaf, it satisfies descent. 
As constructibility is preserved by pullback, we have inclusions of full subcategories
\[
\D_\cons(U,\Lambda) \subset \tot\big(\D_\cons(U_\bullet,\Lambda)\big) \subset \D(U,\Lambda).
\]
The totalization is the full subcategory of objects $M\in \D(U,\Lambda)$ such that $M|_{U_0}$ is constructible.
Hence, $M\in \D_\cons(U,\Lambda)$ by \thref{lisse.constructible.local.property.lemm}.
\xpf

\coro
\thlabel{lisse.evaluation.contractible.lemm}
Let $X$ be a scheme, and let $M\in \D(X,\Lambda)$. 
Then $M$ is lisse if and only if, for every map of w-contractible affines $V\r U$ in $X_\proet$, $\RG(U, M)$ is a perfect complex (over $\Gamma(U, \Lambda)$) and if the natural map 
\[
\RG(U,M)\t_{\Gamma(U,\Lambda)}\Gamma(V,\Lambda) \lr \RG(V,M)
\]
is an equivalence.
\xcoro
\pf
Combine \thref{dualizable.objects.contractible.lem}\refit{dualizable.objects.contractible.lem.3} with \thref{lisse.constructible.hyperdescent.coro}.
\xpf

\subsection{Boundedness properties}

In order to compare our definition of, say, 
constructible $\bar\bbQ_\ell$-sheaves to the classical one in terms of $E$-sheaves for finite field extensions $E\supset \bbQ_\ell$, it is necessary to control hyperdescent not only for $\D_\lis$, but for appropriate colimits of such \ii-categories.
To do this, we filter the \ii-categories of lisse sheaves according to the amplitude of objects:

\defi
\thlabel{bounded.complexes.nota}
For an integer $n\geq 0$, we write $\D_\lis^{\{-n,n\}}(X,\Lambda)$ for the full subcategory of $\D_\lis(X,\Lambda)$ of objects $M$ such that $M$ and its dual $M^\vee$ lie in degrees $[-n,n]$ with respect to the t-structure on $\D(X, \Lambda)$.
\xdefi

The purpose of introducing this  subcategory is to have a $(2n+1)$-category: $M\cong \tau^{\leq n}\tau^{\geq -n}M$ for each such object.
For example, the category $\D_\lis^{\{0,0\}}(X,\Lambda)$ is the full subcategory of dualizable objects in $\D(X,\Lambda)^\heartsuit$, that is, those $\Lambda$-sheaves $M$ that are locally on $X_\proet$ isomorphic to $\underline N\t_{\underline{\Lambda_*}}\Lambda_X$ for some finite projective $\Gamma(X,\Lambda)$-module $N$.

\lemm 
\thlabel{bounded.complexes.lemm}
Assume that $X$ is qcqs.
Then 
\[
\D_\lis(X,\Lambda)=\bigcup_{n\geq 0}\D_\lis^{\{-n,n\}}(X,\Lambda)
\]
as full subcategories of $\D_\lis(X,\Lambda)$.
\xlemm
\pf
The condition of being in the subcategory $\D_\lis^{\{-n,n\}}(X,\Lambda)$ can be checked pro\'etale locally: the restriction functors are monoidal, conservative and preserve the t-structure.
So we may assume that $X$ is a w-contractible qcqs scheme. 
Then, under the equivalence $\D_\lis(X, \Lambda) \cong \Perf_{\Gamma(X, \Lambda)}$ (see \thref{dualizable.objects.contractible.lem} \refit{dualizable.objects.contractible.lem.3}), an object lies in the subcategory $\D_\lis^{\{-n,n\}}(X, \Lambda)$ if and only if it is represented by a bounded complex of finitely generated projective $\Gamma(X, \Lambda)$-modules that is concentrated in degrees $[-n,n]$.
Hence, the lemma follows from the corresponding filtration $\Perf_{\Gamma(X, \Lambda)}=\bigcup_{n\geq 0}\Perf_{\Gamma(X, \Lambda)}^{\{-n,n\}}$ on perfect modules. 
\xpf

\coro
\thlabel{constructible.bounded.coro}
Every constructible $\Lambda$-sheaf on a qcqs scheme is bounded. 
\xcoro
\pf
By an induction on the finite number of strata witnessing the constructibility, using the conservativity and t-exactness of the pair of functors $(j^*,i^*)$ in the notation of \eqref{adjoints.closed.open.immersion.eq}, one reduces to the case of lisse sheaves. 
So we are done by \thref{bounded.complexes.lemm}.
\xpf

\lemm
\thlabel{bounded.lisse.objects.lemm}
Let $n\geq 0$ be an integer.
\begin{enumerate}
\item 
\label{item--bounded.lisse.objects.lemm.1}
For any map $f\co Y\r X$ of schemes the pullback functor $f^*$ restricts to a functor
\[
f^*\co \D_\lis^{\{-n,n\}}(X,\Lambda)\r \D_\lis^{\{-n,n\}}(Y,\Lambda).
\]

\item
\label{item--bounded.lisse.objects.lemm.2}
For any map of condensed rings $\Lambda\r \Lambda'$ the base change functor $(\str)\t_{\Lambda_X}\Lambda'_X$ restricts to a functor
\[
(\str)\t_{\Lambda_X}\Lambda'_X\co \D_\lis^{\{-n,n\}}(X,\Lambda)\r \D_\lis^{\{-n,n\}}(X,\Lambda').
\]

\item
\label{item--bounded.lisse.objects.lemm.3}
The functor $X\mapsto \D_\lis^{\{-n,n\}}(X,\Lambda)$ satisfies hyperdescent on $X_\proet$.
\end{enumerate}
\xlemm
\pf
Part (1) is clear since $f^*$ is t-exact and, being monoidal, preserves duals.
For (2) we use that $\t$ is right t-exact in general. 
On the other hand, 
\[
M \mapsto M \t_{\Lambda_X} \Lambda'_X = (M^{\vee})^\vee \t_{\Lambda_X} \Lambda'_X = \IHom_{\Lambda_X}(M^\vee, \Lambda'_X) 
\]
is also left t-exact.
Part (3) is immediate from \thref{lisse.constructible.hyperdescent.coro} and (1), using that the condition of lying in the subcategory $\D_\lis^{\{-n,n\}}(X,\Lambda)$ can be checked pro\'etale locally. 
\xpf

\subsection{Local constancy of lisse sheaves}
\label{sect--lisse.locally.connected.sect}

Recall from \thref{perfect.lc} that a sheaf $M \in \D(X, \Lambda)$ is called pro\'etale-locally perfect-constant if $M$ is locally on $X_\proet$ isomorphic to $\underline{N}\t_{\underline{\Lambda_*}}\Lambda_X$ for some $N\in \Perf_{\Lambda_*}$, where $\Lambda_*=\Gamma(*,\Lambda)$ is the underlying ring.

\theo
\thlabel{locally.constant.prop}
Let $\Lambda$ be a condensed ring.
Let $X$ be a scheme that has locally a finite number of irreducible components.  
Then $M\in \D(X,\Lambda)$ is lisse if and only if $M$ is pro\'etale-locally perfect-constant.
\xtheo
\pf
Let $M$ be lisse (the other direction is clear).
After a Zariski localization, we reduce to the case where $X$ is affine and connected with finitely many irreducible components.
As any two points of $X$ can be joined by a finite zig-zag of specializations, the pullback of $M$ to any geometric point is perfect-constant (\thref{dualizable.objects.contractible.lem}) with the same value $N\in \Perf_{\Lambda_*}$. 
Let $U\in X_\proet$ be any w-contractible affine cover.
We claim that there exists an isomorphism $M|_U\simeq \underline{N}\t_{\underline{\Lambda_*}}\Lambda_U$, implying the theorem.

First, assume that $X$ is irreducible, and fix a geometric generic point $\eta \r X$. 
Let $U_\eta:=U\x_X\eta$ be the base change, and consider the commutative diagram of sites
\eqn
\label{local.constancy.diagram.eq}
\xymatrix{
\eta_\proet \ar[d]  & (U_\eta)_\proet \ar[d]\ar[r]_{\cong}^{\pi_{U_\eta}}\ar[l] & (\pi_0U_\eta)_\proet \ar@{->>}[d]\\
X_\proet & U_\proet \ar[r]^{\pi_U} \ar[l]& (\pi_0U)_\proet \ar@{.>}@/^/[u].
}
\xeqn
Here $\pi_{U_\eta}$ is an equivalence because $\eta$ is a geometric point and $U_\eta\r \eta$ is proétale. 
Further, the map $\pi_0 U_\eta\r \pi_0U$ is surjective and admits a splitting:
the map on topological spaces $|U_\eta|\r |U|\x_{|X|}|\eta|$ is surjective by \StP{03H4}, and hence induces a surjection on connected components.
So we need to see that the image of every connected component of $|U|$ under the map $|U|\r |X|$ contains the unique generic point of the irreducible space $|X|$. 
This is true because the map $|U|\r |X|$ is generalizing and connected components are closed under generalizations.
Now as $U$ is w-contractible affine, so $\pi_0U$ is extremally disconnected profinite, there exists a section to the surjection $\pi_0 U_\eta\r \pi_0U$.

To finish the argument in the irreducible case, we apply $\D_\lis(\str,\Lambda)$ to the diagram \eqref{local.constancy.diagram.eq} and observe that $\pi_U$ induces an equivalence by \thref{connected.component.pullback} \refit{connected.component.pullback.3}.
More concretely, any isomorphism $M|_\eta\simeq \underline{N}\t_{\underline{\Lambda_*}}\Lambda_\eta$ induces an isomorphism over $U_\eta$, and hence over $U$ by using the section.

Next, if $X=\cup X_i$ is the union of finitely many irreducible components, then we denote $U_i=X_i\x_X U$. 
It follows from the irreducible case that there exist isomorphisms $\varphi_i\co M|_{U_i}\simeq \underline{N}\t_{\underline{\Lambda_*}}\Lambda_{U_i}$.
As $U$ is w-contractible affine, so are the closed subschemes $U_i\subset U$ by \cite[Lemma 2.2.15]{BhattScholze:ProEtale}. 
We denote by $\tilde U_i \r U$ their Henselizations.
Using \thref{neighborhood.extension.coro} \refit{neighborhood.extension.coro.2}, the isomorphisms $\varphi_i$ uniquely extend to $\tilde \varphi_i\co M|_{\tilde U_i}\simeq \underline{N}\t_{\underline{\Lambda_*}}\Lambda_{\tilde U_i}$.
As there are only finitely many irreducible components, the disjoint union $\sqcup \tilde U_i\r U$ is a cover in $X_\proet$, and hence admits a section because $U$ is w-contractible affine. 
Therefore, we can pullback the isomorphism $\sqcup \tilde \varphi_i$ along the section to obtain an isomorphism $M|_U\simeq \underline{N}\t_{\underline{\Lambda_*}}\Lambda_{ U}$ as desired. 
\xpf


\section{Change of coefficients} 
\label{sect--coefficients.sect}
We show that the \ii-category of lisse and constructible sheaves behaves well under certain sequential limits and filtered colimits in the condensed coefficients $\Lambda$. 
Throughout, let $X$ be a scheme. 

\subsection{Sequential limits}
\label{sect--coefficients.limits.sect}
In this section, let $\Lambda=\lim_{i\geq 1} \Lambda_i$ be a sequential limit of condensed rings such that all transition maps $\Lambda_{i+1}\r \Lambda_i$ are surjective with locally nilpotent kernel.
The last condition means that, for all profinite sets $S$, all elements of the kernel of $\Gamma(S,\Lambda_{i+1})\r \Gamma(S,\Lambda_{i})$ are nilpotent. 
We note that $\Lambda_X$ identifies via the natural map with
\[
\Lambda_X\stackrel\cong\lr  \lim \Lambda_{i,X}\cong \on{Rlim}\Lambda_{i,X},
\]
where we use that sequential limits of surjections are exact in a replete topos, see \cite[Proposition 3.1.10]{BhattScholze:ProEtale}.
In the following all limits will be derived unless mentioned otherwise.
Also recall the generalities about limits of stable (idempotent complete) \ii-categories from \refsect{limits.filtered.colimits}.

\prop
\thlabel{limit.coefficients.lemm}
The following natural functors are equivalences: 
\begin{align*}
\D_\lis(X,\Lambda)\stackrel\cong \lr & \ \lim \D_\lis(X,\Lambda_i) \text{ and } \\
\D_\cons(X,\Lambda)\stackrel\cong \lr & \ \lim \D_\cons(X,\Lambda_i).
\end{align*} 
Both limits are formed using $(\str) \t_{\Lambda_{j}} \Lambda_i$ for $j\geq i$.
An inverse functor is given by $\{M_i\} \mapsto \lim M_i$.
\xprop
\pf
We begin with the \ii-categories of lisse sheaves. 
Both functors $X\mapsto \D_\lis(X,\Lambda), \lim \D_\lis(X,\Lambda_i)$ are hypersheaves on $X_\proet$ by \thref{lisse.constructible.hyperdescent.coro}.
So we reduce to the case where $X$ is w-contractible and affine. 
Using \thref{dualizable.objects.contractible.lem} \refit{dualizable.objects.contractible.lem.3}, we get a commutative (up to equivalence) diagram:
\[
\xymatrix{
\D_\lis(X,\Lambda) \ar[r] \ar[d]^{\cong}& \lim \D_\lis(X,\Lambda_i) \ar[d]^{\cong}
\\
\Perf_{\Gamma(X,\Lambda)} \ar[r] & \lim \Perf_{\Gamma(X,\Lambda_i)}
}
\] 
Since $X$ is w-contractible and affine, all transition maps $\Gamma(X,\Lambda_{i+1})\r \Gamma(X,\Lambda_i)$ are surjective with locally nilpotent kernel.
Thus, the lower horizontal functor is an equivalence by \thref{perfect.complexes.colimit.lemm} \refit{perfect.complexes.colimit.lemm.3}.

As for constructible sheaves, we claim that the map 
\[
M\stackrel\cong \lr \lim \left(M\t_{\Lambda_X}\Lambda_{i,X}\right)
\]
is an equivalence for any $M\in \D_\cons(X,\Lambda)$: the functor $\iota^*$ commutes with limits for any constructible locally closed immersion $\iota\co Z \r X$ by \cite[Corollary 6.1.5]{BhattScholze:ProEtale}. 
Using the localization sequences \eqref{adjoints.closed.open.immersion.eq} we may assume that $M$ is lisse where our claim is already proven.
This formally implies the full faithfulness: for any $N\in \D_\cons(X,\Lambda)$, $M\in \D(X,\Lambda)$,
\[
\Hom_{\Lambda_X}(M,N)\cong \lim \Hom_{\Lambda_X}(M,N\t_{\Lambda_X}\Lambda_{i,X})\cong \lim \Hom_{\Lambda_{i,X}}(M\t_{\Lambda_X}\Lambda_{i,X},N\t_{\Lambda_X}\Lambda_{i,X}).
\]
For essential surjectivity, it is enough to show that, for any $\{M_i\}\in \lim \D_\cons(X,\Lambda_i)$ and fixed $j\geq 1$, the limit $\lim M_i\in \D(X_\proet,\Lambda)$ is constructible and the natural map
\begin{equation}\label{reduction.eq}
(\lim M_i)\t_{\Lambda_X}\Lambda_{j,X}\lr M_j
\end{equation}
is an equivalence.
As before, we reduce to the case where $M_j$ is lisse.
We claim that $\lim M_i$ is lisse as well.
Evaluating at any w-contractible affine $U\in X_\proet$ gives 
\begin{equation}\label{reduction.eq.2}
(\lim \RG(U,M_i))\t_{\Gamma(U,\Lambda)}\Gamma(U,\Lambda_{j})\lr \RG(U,M_j)
\end{equation}
Since $\RG(U,M_j)$, and hence $\RG(U,M_1)$, is perfect by \thref{dualizable.objects.contractible.lem} \refit{dualizable.objects.contractible.lem.3}, this map is an equivalence and $\lim \RG(U,M_i)=\RG(U,\lim M_i)$ is perfect by \StP{0CQG}.
Since $U$ was arbitrary, \eqref{reduction.eq} is an equivalence. 
To see that $\lim M_i$ is lisse, it is remains (\thref{lisse.evaluation.contractible.lemm}) to show that the natural map
\[
(\RG(U,\lim M_i))\t_{\Gamma(U,\Lambda)}\Gamma(V,\Lambda)\lr \RG(V,\lim M_i)
\]
is an equivalence for any map of w-contractible affines $V\r U$ in $X_\proet$.
Using the conservativity of $(\str)\t_{\Gamma(V,\Lambda)} \Gamma(V,\Lambda_j)$ on perfect complexes proven in \thref{perfect.complexes.colimit.lemm} \refit{perfect.complexes.colimit.lemm.3}, this follows from \eqref{reduction.eq.2}.
\xpf

\subsection{Filtered colimits}
\label{sect--coefficients.colimits.sect}
In this section, we assume that $\Lambda=\colim \Lambda_i$ is a filtered colimit of condensed rings $\Lambda_i$.
Recall the generalities about filtered colimits of stable (idempotent complete) \ii-categories from \refsect{limits.filtered.colimits}.

\prop
\thlabel{colimit.coefficients.lemm}
If $X$ is qcqs, then the following natural functors are equivalences:
\begin{align*}
\colim \D_\lis(X,\Lambda_i)\stackrel \cong \lr & \ \D_\lis(X,\Lambda) \text{ and } \\
\colim \D_\cons(X,\Lambda_i)\stackrel \cong \lr & \ \D_\cons(X,\Lambda).
\end{align*}
Both filtered colimits are formed using $(\str) \t_{\Lambda_i} \Lambda_j$ for $j \geq i$.
\xprop

\pf
For lisse sheaves, it suffices by \thref{bounded.complexes.lemm} to show that the functor
\[
\colim_i \D_\lis^{\{-n,n\}}(X,\Lambda_i)\lr \D_\lis^{\{-n,n\}}(X,\Lambda)
\]
is an equivalence for any fixed $n\geq 0$. 
Both sides satisfy hyperdescent on $X_\proet$ (\thref{bounded.lisse.objects.lemm} \refit{bounded.lisse.objects.lemm.3} using \thref{filtered.colimits.limits.lemm}), so we may assume that $X$ is w-contractible qcqs.
In this case we have $\D_\lis^{\{-n,n\}}(X,\Lambda)\cong \Perf_{\Gamma(X,\Lambda)}^{\{-n,n\}}$ by \thref{dualizable.objects.contractible.lem} \refit{dualizable.objects.contractible.lem.3}, see also the proof of \thref{bounded.complexes.lemm}.
Since $X$ is qcqs, we have a presentation $\Gamma(X,\Lambda)=\colim\Gamma(X,\Lambda_i)$ as a filtered colimit of rings. 
We conclude using \thref{perfect.complexes.colimit.lemm} \refit{perfect.complexes.colimit.lemm.1}.

As for constructible sheaves we note that for any constructible locally closed immersion $\iota: Z \r X$ and $M\in \D(Z_\proet,\Lambda)$, $N\in \D(X_\proet,\Lambda)$ we have 
\eqn
\label{tensor.extension.by.zero.eq}
\iota_!(M\t_{\Lambda_Z} \iota^*N)\cong \iota_!M\t_{\Lambda_{X}} N 
\xeqn
by \cite[Lemma 6.2.3 (3)]{BhattScholze:ProEtale}.
Applying this with $N=\Lambda_i$ and using standard arguments involving the fiber sequence $j_!j^*\r \id \r i_*i^*$ in the notation of \eqref{adjoints.closed.open.immersion.eq} the essential surjectivity follows from the case of lisse sheaves.
For full faithfulness, it suffices to show (after using standard $\t$-$\Hom$-adjunctions) that, for any $M,N\in \D_\cons(X,\Lambda_i)$, the natural map
\[
\colim_{j\geq i}\Hom_{\Lambda_{X,i}}\left(M,N\t_{\Lambda_{X,i}}\Lambda_{X,j}\right)  \lr \Hom_{\Lambda_{X,i}}\left(M,\colim_{j\geq i}N\t_{\Lambda_{X,i}}\Lambda_{X,j}\right)
\]
is an equivalence. 
By \thref{colimit.commutation.lemm} below, it is enough to show that there exists an integer $n\geq 0$ such that $N\t_{\Lambda_{i,X}}\Lambda_{j,X}\in \D^{\geq -n}(X,\Lambda_i)$ for all $j\geq i$. 
Using that $X$ is qcqs, we can perform an induction on the number of strata of a stratification witnessing the constructibility of $N$.
Applying \eqref{tensor.extension.by.zero.eq} to $\iota_! \iota^*N \t_{\Lambda_{i, X}} \Lambda_{j, X}$ and using the t-exactness of $\iota_!$ (\thref{t.exactness.lemm}), we may then assume that $N$ is lisse.
Then it is lies in $\D_\lis^{\{-n,n\}}(X,\Lambda_i)$ for $n \gg 0$, so that $N\t_{\Lambda_{i,X}}\Lambda_{j,X}$ lies in the same subcategory as well, see \thref{bounded.lisse.objects.lemm} \refit{bounded.lisse.objects.lemm.2}.
\xpf

In the proof, we used the following general lemma. 
An analogous result for étale sheaves is proven in \cite[Lemma 6.3.14]{BhattScholze:ProEtale}:

\lemm
\thlabel{colimit.commutation.lemm}
For any fixed integer $n\in \bbZ$, the following functors commute with filtered colimits with terms in $\D^{\geq n}(X,\Lambda)$:
\begin{enumerate}

\item 
\label{item--colimit.commutation.lemm.1}
$f_*\co \D(X,\Lambda)\r \D(Y,\Lambda)$ for any qcqs map $f\co X\r Y$;

\item 
\label{item--colimit.commutation.lemm.2}
$\underline\Hom_{\Lambda_X}(M,\str)\co \D(X,\Lambda)\r \D(X,\Lambda)$ for any qcqs scheme $X$ and $M\in \D_\cons(X,\Lambda)$.
\end{enumerate}
In particular, under the conditions in \refit{colimit.commutation.lemm.2}, the functor 
\[
\Hom_{\Lambda_X}(M,\str)=\RG\big(X,\underline\Hom_{\Lambda_X}(M,\str)\big)\co \D(X,\Lambda)\r \Mod_{\Gamma(X,\Lambda)}
\]
commutes with such colimits as well. 
\xlemm
\pf
For the final assertion, we apply \refit{colimit.commutation.lemm.1} to the map of sites $f=p_X\co X_\proet\r *_\proet$. 
Then $\RG(X,\str)$ is the composition of the functors   
\[
\D(X,\Lambda)\stackrel {f_*} \lr \D(*,\Lambda)\stackrel {\RG(*,\str)} \lr \Mod_{\Gamma(X,\Lambda)},
\]
and hence commutes with filtered colimits with terms in $\D^{\geq n}(X,\Lambda)$ as well. 
So using \refit{colimit.commutation.lemm.2}, we see that $\Hom_{\Lambda_X}(M,\str)$ commutes with such filtered colimits as well.  
Here we use that any constructible sheaf on a qcqs scheme is bounded (see \thref{constructible.bounded.coro}), so that the functor $\underline\Hom_{\Lambda_X}(M,\str)$ maps $\D^{\geq n}(X,\Lambda)$ into $\D^{\geq m}(X,\Lambda)$ for some $m\leq n$.

Statement \refit{colimit.commutation.lemm.1} is an instance of \cite[Corollary~3.10.5]{Barwick.Glasman.Haine.exodromy}.
We include an argument for the convenience of the reader.
Let $N=\colim N_j$ be a filtered colimit of some sheaves $N_j\in \D^{\geq n}(X,\Lambda)$.
It is enough to show that the natural map
\eqn
\label{cohomology.sheaves.colimit.eq}
\colim_j\H^p\circ f_*(N_j) \r \H^p\circ f_*(N)
\xeqn
is an equivalence in $\D(Y,\Lambda)^\heartsuit$ for any $p\in \bbZ$, $\H^p:=\tau^{\leq p}\circ \tau^{\geq p}$.
As filtered colimits are t-exact we can write $N=\colim N_j=\colim_{m,j} \tau^{\leq m}N_j$. 
By left exactness of $f_*$ only the terms $\tau^{\leq p}N_j, \tau^{\leq p}N$ contribute to $\H^p\circ f_*$ and we may assume $N_j, N\in \D^{[n,p]}(X,\Lambda)$.
An induction on the length $p-n$ reduces us further to the case where $N_j, N$ are in a single t-degree.
So after possibly renumbering we may assume $N_j, N\in \D(X,\Lambda)^\heartsuit$ embedded in degree $0$. 
Evaluating  \eqref{cohomology.sheaves.colimit.eq} at any $V\in Y_\proet$ w-contractible affine it is enough to show that
\[
\colim_j \H^p(X\x_YV,N_j)=\colim_j \Gamma(V,\H^p\circ f_* (N_j))\r \Gamma(V,\H^p\circ f_* (N))=\H^p(X\x_YV,N)
\]
is an isomorphism.
By our assumption on $f$, the base change $X\x_YV$ is qcqs. 
It remains to show that for any qcqs scheme $X$ the cohomology functor $\H^p(X,\str)\co \D(X,\Lambda)^\heartsuit\r \Mod_{\Gamma(X,\Lambda)}^\heartsuit$ commutes with filtered colimits. 
Choosing a hypercover $U_\bullet\r X$ in $X_\proet$ by w-contractible affine schemes, this can be computed as the $p$-th cohomology of the complex
\[
\dots\r 0\r \Gamma(U_0,\str) \r \Gamma(U_1,\str)\r \Gamma(U_2,\str)\r\dots.
\]
As each $\Gamma(U_i,\str)$, $i\geq 0$ commutes with filtered colimits so does $\H^p(X,\str)$.

For \refit{colimit.commutation.lemm.2}, let $X$ be qcqs and $M\in \D_\cons(X,\Lambda)$.
We claim that the functor $\underline\Hom_{\Lambda_X}(M,\str)$ commutes with filtered colimits with terms in $\D^{\geq 0}(X_\proet,\Lambda)$.
If $M$ is lisse (=dualizable), then $\underline\Hom_{\Lambda_X}(M,\str)=(\str)\t_{\Lambda_X}M^\vee$ commutes with all colimits. 
In general, by an induction on the finite number of strata in $X$ witnessing the constructibility of $M$, we reduce to the case $M=\iota_!\iota^*M$ where $\iota^*M$ is lisse for some constructible locally closed immersion $\iota\co Z\hr X$. 
Using standard adjunctions we compute
\[
\underline\Hom_{\Lambda_X}\left(\iota_!\iota^*M,\str\right)=\iota_*\underline\Hom_{\Lambda_Z}\left(\iota^*M,\iota^!(\str)\right)=\iota_*\left(\iota^!(\str)\t_{\Lambda_Z} (\iota^*M)^\vee\right).
\]
Note that $\iota^!$ is left t-exact (as the right adjoint of the t-exact functor $\iota_!$), so preserves the subcategory $\D^{\geq 0}(X,\Lambda)$. 
In light of \refit{colimit.commutation.lemm.1} applied to $\iota_*$, it remains to show that $\iota^!$ commutes with the desired colimits.
If $\iota$ is an open immersion, then the t-exact functor $\iota^!=\iota^*$ commutes with all colimits. 
We reduce to the case where $\iota=i\co Z\hr X$ is a constructible closed immersion with open complement $j\co U\hr X$. 
Then the fiber sequence $i_*i^!\r \id \r j_*j^*$ and \refit{colimit.commutation.lemm.1} applied to $j_*$ shows that $i^!$ commutes with the desired colimits as well. 
\xpf

\rema
\thlabel{cohomological.dimension.warn}
The condition that the filtered colimit is formed using objects in $\D^{\ge n}(X, \Lambda)$ can not in general be dropped in \thref{colimit.commutation.lemm} (see however \thref{all.colimit.commutation.lemm} for a positive result in this direction):
Assume that $\RG(X,\Lambda)$ is concentrated in infinitely many degrees. 
For example, $\RG(\Spec(\bbR),\bbZ/2)$ computes the group cohomology of $\Gal(\bbC/\bbR)=\bbZ/2$ on the trivial module $\bbZ/2$ which is equal to $\bbZ/2$ in all even positive degrees.
Since $\D(X,\Lambda)$ is left-complete \cite[Proposition 3.3.3]{BhattScholze:ProEtale}, the natural map $\oplus_{n\geq 0}\Lambda[n]\to \prod_{n\geq 0}\Lambda[n]$ is an equivalence.
If $\H^0(X,\str)$ commuted with infinite direct sums, we would obtain a contradiction:
\[
\bigoplus_{n\geq 0}\H^0\left(X,\Lambda[n]\right)=\H^0\left(X,\bigoplus_{n\geq 0}\Lambda[n]\right)=\H^0\left(X,\prod_{n\geq 0}\Lambda[n]\right)=\prod_{n\geq 0}\H^0\left(X,\Lambda[n]\right).
\]
\xrema

\subsection{Localizations}
\label{sect--coefficients.localizations.sect}
In this section, let $\Lambda$ be a condensed ring and $T\subset \Gamma(*,\Lambda)$ a multiplicatively closed subset. 
Then the localization $T^{-1}\Lambda=\colim_{t\in T}\Lambda$ viewed as a filtered colimit of sheaves defines a condensed ring.
Since $*$-pullbacks commute with colimits, we have $T^{-1}\Lambda_X=\colim_{t\in T}\Lambda_X$.
Its values on any qcqs $U\in X_\proet$ are computed as
\[
\Gamma(U, T^{-1}\Lambda)=\colim_{t\in T}\Gamma(U, \Lambda) =T^{-1}\Gamma(U,\Lambda).
\]
The second equality is clear, and the first equality is an instance of \eqref{cohomology.sheaves.colimit.eq}.
Let $\Lambda_*:=\Gamma(*,\Lambda)$, and denote by $T^{-1}\Lambda_*$ its localization.

\prop
\thlabel{localization.lisse.cons.lemm}
If $X$ is qcqs, the following functors induced by $M\t_{\Lambda_*} T^{-1}\Lambda_*\mapsto M\t_{\Lambda_X} T^{-1}\Lambda_X$ are fully faithful: 
\begin{align*}
\D_\lis(X,\Lambda)\t_{\Perf_{\Lambda_*}} \Perf_{T^{-1}\Lambda_*} \r &  \ \D_\lis\left(X,T^{-1}\Lambda\right)
\text{ and }  \\
\D_\cons(X,\Lambda)\t_{\Perf_{\Lambda_*}} \Perf_{T^{-1}\Lambda_*} \r & \ \D_\cons\left(X,T^{-1}\Lambda\right).
\end{align*}
\xprop
\pf
For $M\in \D_\cons(X,\Lambda)$, we denote $T^{-1}M=T^{-1}\Lambda_X\t_{\Lambda_X}M$. 
Using that $\Perf_{T^{-1}\Lambda_*}$ is generated under finite colimits by $T^{-1}\Lambda_*$, it is enough to show that the natural map
\[
\Hom_{\Lambda_X}(M,N)\t_{\Lambda_*}T^{-1}\Lambda_* \r \Hom_{\Lambda_X}\big(M,T^{-1}N\big)=\Hom_{T^{-1}\Lambda_X}\big(T^{-1}M,T^{-1}N\big)
\]
is an equivalence for any $M,N\in \D_\cons(X,\Lambda)$.
This follows from \thref{colimit.commutation.lemm}.
\xpf

The functor from $T1$-topological abelian to condensed abelian groups does not commute with filtered colimits in general.
However, the following lemma shows, for example, that $\bbQ_\ell=\colim_{\ell\x}\bbZ_\ell$ and that $\bar\bbQ_\ell=\colim_{E/\bbQ_\ell \,\on{finite}} E$ (writing each $E$ as a filtered colimit of $\calO_E$'s) holds as condensed rings:

\lemm
\thlabel{filtered.colimit.Hausdorff.rings}
Let $\Lambda=\colim \Lambda_i$ be a countable filtered colimit of quasi-compact Hausdorff topological abelian groups with injective transition maps.
Then the induced map of condensed abelian groups $\colim_i\Lambda_i\r \Lambda$ is an isomorphism.
\xlemm
\pf
First off, filtered colimits exist in the category of topological abelian groups (or topological rings) and are formed by taking the colimit in the category of abelian groups (or rings) equipped with its colimit topology.
It is enough to show that the map $\colim \Gamma(S,\Lambda_i)\r \Gamma(S,\Lambda)$ is an isomorphism for any profinite set $S$.
Injectivity is clear. 
For surjectivity, we claim that every continuous map $S\r \Lambda$ factors through some $\Lambda_i$.
As injections between quasi-compact Hausdorff spaces are closed embeddings, this follows from \cite[Lemma 4.3.7]{BhattScholze:ProEtale}.
\xpf

\rema
\thlabel{etale.vs.proetale.fundamental.group}
The functor between the \ii-categories of lisse sheaves in \thref{localization.lisse.cons.lemm} is not an equivalence in general as the source category does not satisfy descent, compare with the discussion above \cite[Lemma 7.4.7]{BhattScholze:ProEtale} (see, however, \thref{classical.comparison.lemm} for a positive result for constructible sheaves). 
More precisely, the failure of essential surjectivity accounts for the difference between the \'etale and pro\'etale fundamental group. 
For example, let $X=\bbP^1/0\sim \infty$ be the nodal curve over some algebraically closed field. 
Its pro\'etale fundamental group (with respect to the choice of some geometric point) is $\pi_1^\proet(X)=\bbZ$ equipped with the discrete topology, whereas the \'etale fundamental group is its profinite completion $\pi_1^\et(X)=\widehat\bbZ$.
So the category $\D_\lis^{\{0,0\}}(X,\Ql)$ is the category of (continuous) representations of $\bbZ$ on finite-dimensional $\bbQ_\ell$-vector spaces, whereas the source category corresponds to the strict full subcategory of those representations stabilizing a $\bbZ_\ell$-lattice. 
\xrema


\section{t-Structures}
\label{sect--t.structures.sect}
The definition of lisse and constructible sheaves is well-adapted to the derived setting. 
The natural t-structure on the \ii-category of all sheaves only restricts to a t-structure on the \ii-categories of lisse and constructible sheaves under additional assumptions on the scheme $X$ and the condensed ring of coefficients $\Lambda$. 
We denote by 
\eqn
\label{t.structures.eq}
\D_\lis^{\geq 0}(X,\Lambda):=\D_\lis(X,\Lambda)\cap \D^{\geq 0}(X,\Lambda) \text{ and } \D_\cons^{\geq 0}(X,\Lambda):=\D_\cons(X,\Lambda)\cap \D^{\geq 0}(X,\Lambda)
\xeqn
the full subcategories of $\D_\lis(X,\Lambda)$, respectively $\D_\cons(X,\Lambda)$, and likewise for the subcategories in cohomological degrees $\leq 0$.
Following \cite[Chapter 6, Section 2]{Glaz:CoherentRings}, we say that a (unital, commutative) ring is {\it regular coherent} if every finitely generated ideal is finitely presented (that is, the ring is coherent \StP{05CU}) and has finite projective dimension.

\defi
\thlabel{t.admissible.defi}
A condensed ring $\Lambda$ is called {\it t-admissible} if the underlying ring $\Lambda_*=\Gamma(*,\Lambda)$ is regular coherent and $\Lambda_*\r \Gamma(S,\Lambda)$ is flat for any extremally disconnected profinite set $S$.
\xdefi

We show in \thref{t.structure.discrete} that the flatness condition is automatic for a discrete topological ring $\Lambda$ (viewed as a condensed ring as in \thref{topological.rings.exam}).
Thus, $\Lambda$ is t-admissible if and only if $\Lambda_*$ is regular coherent.
For example, this holds if $\Lambda_*$ is regular Noetherian of finite Krull dimension, see \thref{t.admissible.properties} \refit{t.admissible.properties.4}.  
Further examples of t-admissible condensed rings include all T1-topological rings such that $\Lambda_*$ is semi-hereditary, see \thref{semi.hereditary.t.admissible}.
This covers algebraic field extensions of $E\supset \bbQ_\ell$ for some prime $\ell$ and their rings of integers $\calO_E$, but also more exotic choices such as the real and complex numbers $\bbR$, $\bbC$ with their Euclidean topology and the ring of adeles $\bbA_K$ for some number field $K$, see \thref{examples.t.admissible.rings}.

\theo
\thlabel{t.structure.condensed}
Let $\Lambda$ be a condensed ring.
\begin{enumerate}
\item 
The natural t-structure on $\D(*,\Lambda)$ restricts to a t-structure on $\D_\lis(*,\Lambda)$ if and only if $\Lambda$ is t-admissible.

\item
\label{item--t.structure.condensed.1}
If $\Lambda$ is t-admissible and $X$ has Zariski locally a finite number of irreducible components, the \ii-categories in \eqref{t.structures.eq} define a t-structure on $\D_\lis(X,\Lambda)$.
Its heart $\D_\lis(X,\Lambda)^\heartsuit$ is the full subcategory of $\D(X,\Lambda)^\heartsuit$ of sheaves $M$ that are locally on $X_\proet$ isomorphic to $\underline {N}\otimes_{\underline{\Lambda_*}}\Lambda_X$ for some finitely presented $\Lambda_*$-module $N$.

\item 
\label{item--t.structure.condensed.2}
If $\Lambda$ is t-admissible and every constructible subset in $X$ has Zariski locally finitely many irreducible components, then the \ii-categories in \eqref{t.structures.eq} define a t-structure on $\D_\cons(X,\Lambda)$.
Its heart $\D_\cons(X,\Lambda)^\heartsuit$ is the full subcategory of $\D(X,\Lambda)^\heartsuit$ of sheaves $M$ such that, for any open affine $U\subset X$, there exists a finite subdivision of $U$ into constructible locally closed subschemes $U_i\subset U$ such that $M|_{U_i}\in \D_\lis(X,\Lambda)^\heartsuit$.
\end{enumerate}
\xtheo

\coro\thlabel{locally.constant.cor}
Let $\Lambda$ be a t-admissible condensed ring.
Let $X$ be a qcqs scheme having locally a finite number of irreducible components.  
Then $M\in \D(X,\Lambda)$ is lisse if and only if $M$ is bounded and each cohomology sheaf is locally on $X_\proet$ isomorphic to $\underline {N}\otimes_{\underline{\Lambda_*}}\Lambda_X$ for some finitely presented $\Lambda_*$-module $N$.
\xcoro
\pf
If $M$ is lisse, then $M$ is bounded (as $X$ is qcqs, see \thref{bounded.complexes.lemm}) and each cohomology sheaf $\H^p(M)$, $p\in \bbZ$ is lisse, using the t-admissibility of $\Lambda$ (\thref{t.structure.condensed}).
The converse follows from an easy induction on the length using that $\D_\lis(X,\Lambda)$ is stable.
\xpf

\rema
\thlabel{not.t.structure.rema}
Some finiteness assumption on $X$ is necessary in order to have a t-structure on $\D_\lis(X,\Lambda)$ such that the inclusion into $\D(X,\Lambda)$ is t-exact. 
As a concrete example take $X=\beta \bbN$, $\Lambda=\bbQ_\ell$ and $f\in \Gamma(X,\Lambda)$ as in \thref{not.t.exact.example}.
Let $K$ be the kernel of $f\co \Lambda_X\r \Lambda_X$ formed in $\D(X,\Lambda)^\heartsuit$.
Then $\RG(X,K)=\Gamma(X,K)=0$, but $K\neq 0$ as its stalks at the boundary $\partial X = \beta \bbN\backslash \bbN$ are non-zero.
When combined with the equivalence $\D_\lis(X,\Lambda)\cong \Perf_{\Gamma(X,\Lambda)}$ from \thref{dualizable.objects.contractible.lem} \refit{dualizable.objects.contractible.lem.3}, this shows that $K$ is not lisse. 
In view of \thref{lisse.evaluation.contractible.lemm}, the failure is explained by the lack of the depicted base change property in this corollary.
As a warning, let us point out that $\D_\lis(X,\Lambda)\cong\Perf_{\Gamma(X,\Lambda)}$ inherits the t-structure from $\Mod_{\Gamma(X,\Lambda)}$ because $\Gamma(X,\Lambda)=\Maps_\cont(\beta \bbN, \bbQ_\ell)$ is regular coherent, see \thref{semi.hereditary.sections.combined}. 
However, if one equips $\D_\lis(X,\Lambda)$ with this t-structure, the inclusion into $\D(X,\Lambda)$ will not be t-exact.
\xrema

The proof of \thref{t.structure.condensed} relies on the following key characterization of regular coherent rings. 
We first provide a well-known auxiliary lemma:

\lemm
\thlabel{coherent.ring.Mod.fp.abelian}
A ring $\Lambda$ is coherent if and only if the subcategory $\Mod_\Lambda^{\heartsuit,\fp} \subset \Mod_\Lambda^{\heartsuit}$ (of the abelian category of $\Lambda$-modules) spanned by the finitely presented $\Lambda$-modules is abelian.
\xlemm

\pf
If $\Lambda$ is coherent, then $\Mod_\Lambda^{\heartsuit,\fp}$ is abelian by \StP{05CW}. 
Conversely, assume that $\Mod_\Lambda^{\heartsuit,\fp}$ is abelian.
If $I\subset \Lambda$ is a finitely generated ideal, then $\Lambda\r \Lambda/I$ is a map of finitely presented $\Lambda$-modules, and hence $I=\ker(\Lambda\r \Lambda/I)$ is finitely presented as well. 
Here we used that the inclusion $\Mod_\Lambda^{\heartsuit,\fp}\subset \Mod_\Lambda^{\heartsuit}$ is left-exact.
This shows that $\Lambda$ is coherent. 
\xpf

\prop
\thlabel{restriction.t.structures}
For any ring $\Lambda$, the following are equivalent: 
\begin{enumerate}
\item 
\label{item--restriction.t.structures.1}
The natural t-structure on $\Mod_\Lambda$ restricts to a t-structure on $\Perf_\Lambda$. 
\item 
\label{item--restriction.t.structures.2}
The ring $\Lambda$ is coherent and every finitely presented $\Lambda$-module has finite Tor dimension. 

\item
\label{item--restriction.t.structures.3}
The ring $\Lambda$ is coherent and every finitely generated ideal has finite Tor dimension. 

\item
\label{item--restriction.t.structures.4}
The ring $\Lambda$ is regular coherent.

\end{enumerate}
\xprop
\pf
We show that \refit{restriction.t.structures.1} implies \refit{restriction.t.structures.2}. By the assumption \refit{restriction.t.structures.1} there are inclusions
\eqn
\label{restriction.t.structures.eq}
\Mod_\Lambda^{\heartsuit,\fp}\subset \Perf_\Lambda^\heartsuit \subset \big(\Mod_\Lambda^\heartsuit\big)^\omega
\xeqn
of full subcategories of $\Mod_\Lambda^\heartsuit$ of finitely presented $\Lambda$-modules at the left and of compact objects at the right:
The first inclusion means that every finitely presented $\Lambda$-module is of the form $\H^0(M)$ for some $M\in \Perf_\Lambda$ which is clear. 
For the second inclusion, we note that any perfect complex is compact in $\Mod_\Lambda$ \StP{07LT} and that the inclusion $\Mod_\Lambda^\heartsuit\subset \Mod_\Lambda$ is full and preserves filtered colimits.
It is well-known that the categories at the left and at the right in \eqref{restriction.t.structures.eq} agree.
Thus, both inclusions are equalities. 
Being the heart of a t-structure, $\Mod_\Lambda^{\heartsuit,\fp}$ is abelian, so that $\Lambda$ is coherent (\thref{coherent.ring.Mod.fp.abelian}).
The inclusion $M\in \Mod_\Lambda^{\heartsuit,\fp}=\Perf_\Lambda^\heartsuit\subset \Perf_\Lambda$ shows that every finitely presented $\Lambda$-module $M$ is perfect. 
By \StP{0658}, this is equivalent to $M$ being pseudo-coherent (or, almost perfect) and of finite Tor dimension. 
This implies \refit{restriction.t.structures.2}.

Conversely, assume that \refit{restriction.t.structures.2} holds.
Then $\Mod_\Lambda^{\heartsuit,\fp}$ is abelian (\thref{coherent.ring.Mod.fp.abelian}) so that every $\Lambda$-module of the form $\H^0(M)$, $M\in \Perf_\Lambda$ is finitely presented. 
Also, every finitely presented $\Lambda$-module is pseudo-coherent by \StP{0EWZ} and, hence perfect since it has finite Tor dimension \StP{0658}.
So $\Perf_\Lambda$ is stable under the truncation functors $\tau^{\leq n}$, $\tau^{\geq n}$ for all $n\in \bbZ$.
This implies \refit{restriction.t.structures.1} since the other properties of a t-structure are inherited from $\Mod_\Lambda$.

It is clear that \refit{restriction.t.structures.2} implies \refit{restriction.t.structures.3}.
We now show the converse implication.
Let $M$ be a finitely presented $\Lambda$-module. 
We need to show that there exists an integer $n>0$ (possibly depending on $M$) such that $\H^p(N\t_\Lambda M)=0$ for all $p>n$ and $N\in \Mod_\Lambda^\heartsuit$.
The argument is similar to the proof of \StP{00HD}:
As $M$ is finitely presented, there is some $m\geq 1$ and an exact sequence $0\r M'\r \Lambda^m\r M \r 0$.
Then $M'$ is finitely presented as well because $\Lambda$ is coherent (so $\Mod_\Lambda^{\heartsuit, \fp}$ is abelian). 
We reduce to the case where $M\subset R^m$ is a submodule.
If $m=1$, then $M$ is a finitely generated ideal and we are done.
If $m\geq 2$, then there is an exact sequence 
\[
0\r M' \r M \r M''\r 0,
\]
where $M'=M\cap (R\oplus 0^{m-1}) \subset R$ and $M''\subset R^{m-1}$ are submodules.   
By induction, there are finitely many finitely generated ideals in $R$ whose Tor dimension bound the Tor dimension of $M$.
This implies \refit{restriction.t.structures.2}.

It remains to prove the equivalence of \refit{restriction.t.structures.3} and \refit{restriction.t.structures.4}.
If $\Lambda$ is as in \refit{restriction.t.structures.3}, then every finitely generated ideal admits a finite resolution by finite projective modules, using the equivalent characterization \refit{restriction.t.structures.1}.
Thus, \refit{restriction.t.structures.3} implies \refit{restriction.t.structures.4}.
Conversely, any finite projective resolution is K-flat.
So ideals admitting such resolutions are of finite Tor dimension, proving \refit{restriction.t.structures.4} implies \refit{restriction.t.structures.3}.
\xpf

\pf[Proof of \thref{t.structure.condensed}]
First, assume that $\Lambda$ is t-admissible and that $X$ has locally a finite number of irreducible components, respectively every constructible subset has so. 
We show that the categories $\D_\lis(X,\Lambda)$, respectively $\D_\cons(X,\Lambda)$ are closed under the truncation functors $\tau^{\leq 0}$, and hence inherit the t-structure from $\D(X,\Lambda)$. 
Since restriction commutes with truncation functors, we reduce to the case of $\D_\lis(X,\Lambda)$ with $X$ being affine and connected with finitely many irreducible components. 
So pick $M\in \D_\lis(X,\Lambda)$.
We need to show that $\tau^{\leq 0}M$ is lisse as well. 
For any w-contractible affine cover $U\in X_\proet$, there is an isomorphism $M|_U\simeq \underline N\t_{\underline{\Lambda_*}}\Lambda_U$ for some $N\in \Perf_{\Lambda_*}$, see \thref{locally.constant.prop}.
We compute
\[
(\tau^{\leq 0}M)|_U \cong \tau^{\leq 0}M|_U \simeq \tau^{\leq 0}\left(\underline N\t_{\underline{\Lambda_*}}\Lambda_U\right) \stackrel\cong\lr {\underline{\tau^{\leq 0}N}}\t_{\underline{\Lambda_*}}\Lambda_U,
\]
where the last map is checked to be an isomorphism by evaluating at w-contractible affines $V\in U_\proet$ and using the flatness of $\Lambda_*\r \Gamma(V,\Lambda)=\Gamma(\pi_0V,\Lambda)$.
Note that $\RG(V,\str)$ is t-exact by \thref{dualizable.objects.contractible.lem} \refit{dualizable.objects.contractible.lem.1}, that $\pi_0(V)$ is extremally disconnected \cite[Lemma 2.4.8]{BhattScholze:ProEtale}, and that $\Lambda$ is assumed to be t-admissible.
Further, since $\Lambda_*$ is regular coherent, \thref{restriction.t.structures} shows that $\tau^{\leq 0}N\in \Perf_{\Lambda_*}$.
So $\tau^{\leq 0}M$ is pro\'etale-locally perfect-constant, thus lisse. 
Also, the description of the hearts in \refit{t.structure.condensed.1} and \refit{t.structure.condensed.2} follows immediately from \eqref{restriction.t.structures.eq}.

It remains to show that t-admissibility is necessary in order to have the restricted t-structure on lisse sheaves on the point. 
So assume that the natural t-structure on $\D(*,\Lambda)$ restricts to a t-structure on $\D_\lis(*,\Lambda)$.
In particular, the latter category is closed under truncation in $\D(*,\Lambda)$. 
As $\RG(*,\str)$ is t-exact, we see that $\Perf_{\Lambda_*}$ is closed under truncation in $\Mod_{\Lambda_*}$.
By the equivalent characterization in \thref{restriction.t.structures}, the ring $\Lambda_*$ is regular coherent.
Similarly, using the t-exactness of $\RG(S,\str)$ for any extremally disconnected profinite set $S$, we see that the functor $\Perf_{\Lambda_*}\r \Mod_{\Gamma(S,\Lambda)}$, $N\mapsto N\t_{\Lambda_*}\Gamma(S,\Lambda)$ is t-exact.
We claim that $\on{Tor}^{\Lambda_*}_1(\Gamma(S,\Lambda), \Lambda_*/I)=0$ for all finitely generated ideals $I\subset \Lambda_*$ which implies flatness of $\Lambda_*\r \Gamma(S,\Lambda)$ by \StP{00M5}.
Indeed, as $\Lambda_*$ is regular coherent, we see $\Lambda_*/I\in \Perf_{\Lambda_*}$ when placed in cohomological degree $0$, say. 
By assumption, $(\Lambda_*/I)\t_{\Lambda_*}\Gamma(S,\Lambda)$ is concentrated in degree $0$ so that we get the desired vanishing. 
We conclude that $\Lambda$ is t-admissible which finishes the proof. 
\xpf

We now exhibit examples of t-admissible condensed rings. 
Throughout, we freely use the equivalent characterizations of regular coherent rings in \thref{restriction.t.structures}.
Recall that a \emph{semi-hereditary} ring is one where every finitely generated ideal is projective.

\lemm
\thlabel{t.admissible.properties}
There holds:
\begin{enumerate}

\item
\label{item--t.admissible.properties.4}
Regular Noetherian rings of finite Krull dimension and semi-hereditary rings are regular coherent.

\item 
\label{item--t.admissible.properties.1}
Regular coherent rings are stable under the following operations: localizations, finite products and filtered colimits with flat transition maps.

\item
\label{item--t.admissible.properties.5}
Semi-hereditary rings are stable under the following operations: localizations, arbitrary products and filtered colimits.

\end{enumerate}
\xlemm
\pf
Clearly, semi-hereditary rings are regular coherent: every finitely generated projective module is finitely presented.
The remainder of part \refit{t.admissible.properties.4} follows from \StPd{00OE}{0CXE}.

For \refit{t.admissible.properties.1}, let $\Lambda=\colim \Lambda_i$ be a filtered colimit with flat transition maps.
We observe that every finitely generated ideal $I\subset \Lambda=\colim \Lambda_i$ is of the form $I_{j}\t_{\Lambda_j}\Lambda$ for some finitely generated ideal $I_j\subset \Lambda_j$, using the flatness of the transition maps.
Thus, if $I_j$ is finitely presented and of finite Tor dimension as a $\Lambda_j$-module, so is $I$ as a $\Lambda$-module. 
Similarly, given a finitely generated ideal in a localization $I\subset T^{-1}\Lambda$, there exists a finitely generated ideal $J\subset \Lambda$ such that $I=T^{-1}J$.
So, if $J$ is finitely presented and of finite Tor dimension, so is $I$. 
Finite products are easy and left to the reader.

For \refit{t.admissible.properties.5}, we use \cite[Section 4.2, Proposition]{Olivier:AbsolutelyFlat} for products. 
The remaining arguments are as above. Note that $I_j \t_{\Lambda_j} \Lambda$ there is automatically an ideal in $\Lambda$ since $I_j$ is projective, hence flat, so no flatness of $\Lambda_j \r \Lambda$ is needed. 

\xpf

Now let $\Lambda$ be a condensed ring.
Recall that if $\Lambda$ is associated with a topological ring, we have $\Gamma(S,\Lambda)=\Maps_\cont(S,\Lambda)$ for any profinite set $S$, see \thref{topological.rings.exam}.
For discrete rings, the situation is as nice as possible.

\prop
\thlabel{t.structure.discrete}
For a discrete condensed ring $\Lambda$, the map $\Lambda_* \r \Gamma(S, \Lambda)$ is flat for any profinite set $S$. Thus, $\Lambda$ is t-admissible if and only if $\Lambda_*$ is regular coherent.
\xprop

\pf
Write $S=\lim S_i$ as a cofiltered limit of finite sets.
As $\Lambda$ is discrete, we have 
\[
\Gamma(S,\Lambda)=\colim \Gamma(S_i,\Lambda)=\colim \Lambda_*^{S_i}
\]
which is flat, being a filtered colimit of free, hence flat $\Lambda_*$-modules.
\xpf

\rema
\thlabel{Zl.no.t-structure}
Combining the above proof with \thref{t.admissible.properties} \refit{t.admissible.properties.1} shows that, more generally, $\Gamma(S, \Lambda)$ is regular coherent for any profinite set $S$ provided $\Lambda_*$ is so.

The above shows that for any $\ell > 0$, lisse sheaves with $\Z/\ell^2$-coefficients are not closed under truncation, which bars the attempt to construct the t-structure for $\Zl$-coefficients by a direct limit argument.
We do easily get the t-structure for $\Lambda=\Zl$ using the following method, though.
\xrema


\lemm
\thlabel{semi.hereditary.t.admissible}
Let $\Lambda$ be the condensed ring associated with a T1-topological ring. 
If the underlying ring $\Lambda_*$ is semi-hereditary, then $\Lambda$ is t-admissible. 
\xlemm
\pf
\thref{t.admissible.properties} \refit{t.admissible.properties.4} shows that $\Lambda_*$ is regular coherent.
We claim that $\Lambda_*\r \Gamma(S,\Lambda)$ is flat for any profinite set $S$. 
Since $\Lambda$ is associated with a topological ring, the natural map $\Gamma(S,\Lambda)\r \prod_{s\in S}\Lambda_*$, $f\mapsto (f(s))_{s\in S}$ is injective. 
So $\Gamma(S,\Lambda)$ is a torsionless $\Lambda_*$-module in the sense of Bass \cite[Definition 4.64]{Lam:ModulesAndRings}, and hence flat by a theorem of Chase \cite[Theorem 4.67]{Lam:ModulesAndRings}.
\xpf

Recall that a Prüfer domain is, by one of several equivalent definitions, a semi-hereditary domain. All fields are Prüfer domains, as are the rings of integers $\calO_E$ of algebraic extension fields $E \supset \bbQ_\ell$.

\coro
\thlabel{examples.t.admissible.rings}
The condensed rings associated with the following T1-topological rings are semi-hereditary, hence t-admissible:
\begin{enumerate}

\item
\label{item--examples.t.admissible.rings.1}
All T1-topological Pr\"ufer domains. 

\item
\label{item--examples.t.admissible.rings.2}
The adeles $\bbA_{K}$, the finite adeles $\bbA_{K,{\on{f}}}$ and the profinite completion $\widehat\calO_K$ for any finite field extension $K\supset \bbQ$ with ring of integers $\calO_K$.

\end{enumerate}
\xcoro
\pf
\refit{examples.t.admissible.rings.1} follows directly from \thref{semi.hereditary.t.admissible}.
As localizations and products of semi-hereditary rings are semi-hereditary by \thref{t.admissible.properties} \refit{t.admissible.properties.5}, part \refit{examples.t.admissible.rings.2} follows from the formulas $\widehat\calO_K=\prod_{\mathfrak p}\widehat \calO_{K,\mathfrak p}$ and $\bbA_{K,\on{f}}=\widehat \calO_K[(K\backslash \{0\})^{-1}]$ and $\bbA_{K}\simeq \bbA_{K,{\on{f}}}\x \bbR^{r}\x \bbC^s$ for some $r,s \geq 0$.
\xpf

Finally, let us mention the following result due to Brookshear \cite{Brookshear:ProjectivePrimeIdeals}, De Marco \cite{DeMarco:Projectivity}, Neville \cite{Neville:Coherent} and Vechtomov \cite{Vechtomov:RingsOfFunctions} on the structure of rings of continuous functions. 
Its corollary below is used in \cite{HansenScholze:RelativePerversity} (see \thref{classical.comparison.lemm}) to compute the category of constructible $\bbQ_\ell$-sheaves:

\theo
\thlabel{semi.hereditary.sections.combined}
Let $\Lambda$ be one of the following topological rings:
\begin{enumerate}
\item
\label{item--semi.hereditary.sections.combined.1}
The real numbers $\bbR$ or the complex numbers $\bbC$ with their Euclidean topology

\item
\label{item--semi.hereditary.sections.combined.2}
An algebraic field extension $E\supset \bbQ_\ell$ for some prime $\ell$, or its ring of integers $\calO_E$.

\item
\label{item--semi.hereditary.sections.combined.3}
The adeles $\bbA_{K}$, the finite adeles $\bbA_{K,{\on{f}}}$ or the profinite completion $\widehat\calO_K$ for some finite field extension $K\supset \bbQ$ with ring of integers $\calO_K$.
\end{enumerate}
Then, for any extremally disconnected profinite set $S$, the ring $\Gamma(S,\Lambda)=\Maps_\cont(S,\Lambda)$ is semi-hereditary. 
In particular, it is regular coherent.
\xtheo
\pf
For \refit{semi.hereditary.sections.combined.1}, let $\Lambda$ denote either $\bbR$ or $\bbC$.
If $U\subset S$ is open, then its topological closure $\bar U$ is again open \StP{08YI}. 
In particular, $S$ is a $\Lambda$-basically disconnected Tychonoff (=$T_{3{1\over 2}}$) space or, equivalently, every finitely generated ideal in $\Gamma(S,\Lambda)$ is principal and projective, see \cite[Theorems 8.4.3, 8.4.4]{Glaz:CoherentRings} for $\Lambda=\bbR$ and \cite[Theorem~1]{Vechtomov:Distributive}, \cite[Theorem 13.1]{Vechtomov:RingsOfFunctions} for $\Lambda=\bbC$. 

To prove part \refit{semi.hereditary.sections.combined.2}, the key claim is that for a finite extension $E / \Ql$, the ring $\Gamma(S,\calO_E)$ is semi-hereditary, i.e., every finitely generated ideal is principal and projective.
This implies the claim for the quotient field $E$ and also infinite algebraic extensions by \thref{t.admissible.properties}.

First, we show that every ideal generated by two elements $f,g\in \Gamma(S,\calO_E)$ is principal. 
Let us denote by $|\str|\co E\r \bbR_{\geq 0}$ the normalized valuation so that $\calO_E$ is the subring of elements $x\in E$ with $|x|\leq 1$. 
Let $U=\{s\in S\;|\; |f(s)|>|g(s)|\}$ and $V=\{s\in S\;|\; |f(s)|<|g(s)|\}$ which are open disjoint sets in $S$.
As $S$ is extremally disconnected, their topological closures $\bar U, \bar V$ are clopen (=closed and open) and still disjoint \StP{08YK}.
The characteristic function $e_{\bar U}$ on $\bar U$ defines an idempotent in $\Gamma(S,\calO_E)$.
We claim that the ideal $(f,g)$ is the principal ideal generated by $h:=e_{\bar U}f + (1-e_{\bar U})g\in \Gamma(S,\calO_E)$. 
It is clear that $h\in (f,g)$. 
Conversely, we note that $|f(s)|, |g(s)|\leq |h(s)|$ for all $s\in S$ by construction.
Using the equivalent characterizations in \cite[Theorem 12.2]{Vechtomov:RingsOfFunctions}, we obtain functions $a, b\in \Gamma(S,E)$ with $f=a\cdot h$ and $g=b\cdot h$.
Comparing valuations, we necessarily have $|a(s)|, |b(s)|\leq 1$ for all $s\in S$, that is, $a,b\in \Gamma(S,\calO_E)$. 
So $(f,g)=(h)$ as ideals in $\Gamma(S,\calO_E)$.

It remains to show that every principal ideal is projective. 
More generally, the proof of \cite[Theorem 8.4.4]{Glaz:CoherentRings} shows that this holds for any $T1$-topological ring $\Lambda$ without zero divisors:
For any $f\in \Gamma(S,\Lambda)$, there is a short exact sequence 
\eqn
\label{split.sequence.eq}
0\lr \on{Ann}(f) \lr \Gamma(S,\Lambda)\stackrel{ f\x }\lr (f) \lr 0,
\xeqn
where $\on{Ann}(f)=\{g\;|\;g\cdot f=0\}$ is the annihilator ideal of $f$.
We show that \eqref{split.sequence.eq} splits so that $(f)$ is projective, being a direct summand of $\Gamma(S,\Lambda)$.
Let $U=\{s\in S\;|\; f(s)\neq 0\}$. 
As $\Lambda$ is $T1$, this subset is open and its closure $\bar U\subset S$ is clopen. 
The characteristic function $e_{S\backslash \bar U}$ on $S\backslash \bar U$ defines an idempotent in $\Gamma(S,\Lambda)$.
We claim that $\on{Ann}(f)=(e_{S\backslash \bar U})$ which will imply that \eqref{split.sequence.eq} splits.
Clearly, $e_{S\backslash \bar U}\cdot f=0$ by construction of ${S\backslash \bar U}$. 
Conversely, if $g\cdot f =0 $ for some $g\in \Gamma(S,\Lambda)$, then $g|_U =0$ because $\Lambda$ is without zero divisors. 
As $g^{-1}(0)$ is closed in $S$, we still have $g|_{\bar U}=0$. 
Hence, $g=g\cdot e_{S\backslash \bar U}$, that is, $g\in (S\backslash \bar U)$.

For \refit{semi.hereditary.sections.combined.3}, we apply \thref{t.admissible.properties} which shows that 
\[
\Gamma(S, \widehat \calO_K) = \prod_{\mathfrak p}\Gamma(S, \widehat \calO_{K_{\mathfrak p}})
\]
is semi-hereditary. Here $\mathfrak p$ runs through the places of $K$.
Next, writing $\bbA_{K,\on{f}}=\widehat \calO_K[T^{-1}]$, $T=K\backslash \{0\}$ as a localization and using \thref{filtered.colimit.Hausdorff.rings}, we see that $\Gamma(S,\bbA_{K,\on{f}})=\Gamma(S,\widehat\calO_K)[T^{-1}]$ is semi-hereditary. 
Finally, $\bbA_{K}\simeq \bbA_{K,{\on{f}}}\x \bbR^{r}\x \bbC^s$ for some $r,s \geq 0$ implies the result in this case as well. 
\xpf

\coro
\thlabel{insane.lattice.coro}
Let $S$ be an extremally disconnected profinite set.
Let $E\supset \bbQ_\ell$ and $\bbA_{K,\on{f}}$ be as in \thref{semi.hereditary.sections.combined}.
Then every finite projective module over $\Gamma(S,E)$, respectively over $\Gamma(S,\bbA_{K,\on{f}})$ extends to a finite projective module over $\Gamma(S,\calO_E)$, respectively over $\Gamma(S,\widehat \calO_K)$.
\xcoro
\pf
All rings are semi-hereditary by \thref{semi.hereditary.sections.combined}. 
Every finite projective module $M$ over a semi-hereditary ring $R$ is a direct sum of finitely generated (hence projective) ideals \cite[Theorem~2.29]{Lam:ModulesAndRings}. 
Hence, it is enough to show that every finitely generated ideal extends. 
This is clear from $\Gamma(S,E)=\Gamma(S,\calO_E)[\ell^{-1}]$ and $\Gamma(S,\bbA_{K,\on{f}})=\Gamma(S,\widehat\calO_K)[T^{-1}]$, $T=K\backslash \{0\}$ by multiplying with a common denominator of the continuous functions generating the ideal.
\xpf


\section{Comparison results}
\label{sect--comparison.results}
We compare our definition with classical definitions for discrete rings, adic rings and their localizations. 
The upshot is that we recover \cite[Definition 6.3.1, Definition 6.5.1]{BhattScholze:ProEtale} for discrete and adic rings, and \cite[Definition 6.8.8]{BhattScholze:ProEtale} for algebraic extensions of $\bbQ_\ell$. 
In \refsect{examples.comparison.l.adic}, we give some examples and make the connection to more classical approaches \cite[(1.1)]{Deligne:Weil2}, \cite[Exposé XIII, \S4]{Illusie:TravauxGabber}, \StP{0F4M}, see also \cite[II.5, Appendix A]{KiehlWeissauer}. 

\subsection{Discrete rings}
\label{sect--discrete.coefficients}
In this section, let $\Lambda$ be the condensed ring associated with a {\it discrete} topological ring, also denoted by $\Lambda$, see \thref{topological.rings.exam}. 
Then, for any scheme $X$, the sheaf $\Lambda_X=\underline \Lambda$ is the constant sheaf of rings on $X_\proet$ associated with $\Lambda$.
The morphism onto the étale site $\nu\co X_\proet\r X_\et$ induces a pullback functor
\[
\nu^*\co \D(X_\et,\Lambda) \r \D(X,\Lambda).
\]
Recall from \cite[Remark 6.3.27]{CisinskiDeglise:Etale} that an object $M\in \D(X_\et,\Lambda)$ is dualizable if and only if there exists a covering $\{U_i\r X\}$ in $X_\et$ such that each restriction $M|_{U_i}$ is constant with perfect values.
Let us temporarily denote by 
\[
\D_\lis(X_\et,\Lambda)\subset \D_\cons(X_\et,\Lambda) 
\]
the full subcategories of $\D(X_\et,\Lambda)$ of objects which are lisse (=dualizable in $\D(X_\et,\Lambda)$, by definition), respectively Zariski locally lisse along a finite subdivision into constructible locally closed subschemes. 

\prop
\thlabel{discrete.comparison.prop}
For a discrete ring $\Lambda$, the functor $\nu^*$ induces equivalences
\begin{align*}
\D_\lis(X_\et,\Lambda)\stackrel\cong \lr & \ \D_\lis(X,\Lambda) \text{ and }  \\
\D_\cons(X_\et,\Lambda)\stackrel\cong \lr & \ \D_\cons(X,\Lambda).
\end{align*}
\xprop

\rema
\thlabel{ctf.comparison.rema}
If $X$ is qcqs and $\Lambda$ Noetherian, then $\D_\cons(X_\et,\Lambda)=\D_\ctf(X_\et,\Lambda)$ coincides with the full subcategory of $\D(X_\et,\Lambda)$ of constructible $\Lambda$-sheaves of finite Tor-dimension. 
This can be deduced from the characterization in \StP{03TT}.
\xrema

We need some preparation before giving the proof of \thref{discrete.comparison.prop}.

\lemm
\thlabel{discrete.coefficients.lemm}
Recall that $\Lambda$ is discrete. 
For a profinite set $S=\lim S_i$ there is an equivalence
\[
\colim \Perf_{\Gamma(S_i,\Lambda)}\stackrel \cong \lr  \Perf_{\Gamma(S,\Lambda)}.
\]
Here the transition functors are given by $(\str)\t_{\Gamma(S_i,\Lambda)} \Gamma(S_j,\Lambda)$ for $j\geq i$.
\xlemm
\pf
Any continuous map $S\r \Lambda$ factors through some $S\r S_i$ because $S$ is quasi-compact and $\Lambda$ discrete. 
Hence, $\Gamma(S,\Lambda)=\colim\Gamma(S_i,\Lambda)$ is a filtered colimit. 
So we are done by \thref{perfect.complexes.colimit.lemm} \refit{perfect.complexes.colimit.lemm.1}.
\xpf

For $N\in \Mod_\Lambda$, recall that $\underline N$ denotes the associated constant sheaf of $\Lambda$-modules on $X_\proet$.
If $N\in \Perf_\Lambda$, then $\RG(X,\underline N)\cong\RG(X,\Lambda)\otimes_{\Lambda_*} N$ shows that
\eqn
\label{constant.comparison.eq}
\underline N \cong \left(N\otimes_{\Lambda_*}\Gamma(X,\Lambda)\right)_X\in \D(X,\Lambda).
\xeqn

\coro
\thlabel{locally.constant.lemm}
For a discrete coefficient ring $\Lambda$, an object $M\in \D(X,\Lambda)$ is lisse if and only if there exists a covering $\{U_i\r X\}$ in $X_\proet$ such that each restriction $M|_{U_i}$ is constant with perfect values. 
\xcoro
\pf
Let us assume that $M$ is lisse (the other direction is clear). 
By \thref{lisse.constructible.local.property.lemm}, 
we may assume that $X$ is w-contractible and affine. 
In this case, $P := \pi_0X$ is a profinite set, say $P = \lim P_i$. 
Then \thref{discrete.coefficients.lemm} together with \thref{dualizable.objects.contractible.lem} \refit{dualizable.objects.contractible.lem.3} give equivalences 
\[
\D_\lis(X,\Lambda)\cong \Perf_{\Gamma(X,\Lambda)}=\Perf_{\Gamma(P,\Lambda)}\cong \colim \Perf_{\Gamma(P_i,\Lambda)},
\]
using that $\Gamma(X,\Lambda)=\Gamma(P,\Lambda)$.
In down to earth terms, we find some $N\in \Perf_{\Gamma(P_i,\Lambda)}$ together with an isomorphism
\[
\left(N\t_{\Gamma(P_i,\Lambda)}\Gamma(X,\Lambda)\right)_X\cong M. 
\]
The fibers of the projection $X\r P \r P_i$ induce a finite subdivision of $X$ into clopen subsets.
After Zariski localizing on $X$ we may therefore assume that $P_i=*$.
It follows from \eqref{constant.comparison.eq} that $\underline N\cong M$ is constant.
\xpf

\pf[Proof of \thref{discrete.comparison.prop}]
Since $\nu^*\co \D(X_\et,\Lambda)\r \D(X,\Lambda)$ is monoidal, it induces functors
\eqn
\label{discrete.comparison.prop.proof}
\D_\lis(X_\et,\Lambda)\r \D_\lis(X,\Lambda) \text { and } \D_\cons(X_\et,\Lambda) \r \D_\cons(X,\Lambda),
\xeqn
compatible with étale localization on $X$. 
By étale descent (Zariski descent is enough), we may assume that $X$ is affine. 
The functor $\nu^*$ is fully faithful when restricted to bounded below complexes by \cite[Corollary 5.1.6]{BhattScholze:ProEtale}.
Since $X$ is quasi-compact, any object in $\D_\lis(X_\et,\Lambda)$ and hence in $\D_\cons(X_\et,\Lambda)$ is bounded.
This implies the full faithfulness of the functors in \eqref{discrete.comparison.prop.proof}.

We first show the essential surjectivity for lisse sheaves. 
\thref{locally.constant.lemm} shows that any $M\in \D_\lis(X,\Lambda)$ is perfect-locally constant on $X_\proet$.  
Clearly, any constant sheaf in $\D(X,\Lambda)$ arises as $\nu^*$-pullback from $\D(X_\et,\Lambda)$, that is, it is classical in the terminology of \cite[\S5]{BhattScholze:ProEtale}.
Therefore, locally on $X_\proet$, $M$ is classical.
For sheaves in a single degree, this holds by \cite[Lemma 5.1.4]{BhattScholze:ProEtale}.
The same proof shows this claim for general bounded (below) complexes, see also \cite[above Remark 5.1.7]{BhattScholze:ProEtale}.
Finally, by \cite[Lemma 6.3.13]{BhattScholze:ProEtale} the sheaf $M$ necessarily arises as pullback of a locally on $X_\et$ constant sheaf with perfect values, that is, from an object in $\D_\lis(X_\et,\Lambda)$.

It remains to show the essential surjectivity on constructible sheaves.
If $\iota\co Z\hr X$ is a constructible locally closed subset, the functors $\iota^*, \iota_!$ commute with $\nu^*$, see \cite[Lemma 6.2.3 (4)]{BhattScholze:ProEtale} for $\iota_!$.
Using the full faithfulness of $\nu^*$ on bounded complexes, we reduce by induction on the number of constructible locally closed strata to the case of lisse sheaves.
\xpf

\subsubsection{Adic rings}
\label{sect--adic.coefficients}
In this section, let $\Lambda$ be the condensed ring associated with a Noetherian ring, also denoted by $\Lambda$, complete for the topology defined by an ideal $I\subset \Lambda$, see \thref{topological.rings.exam}.
Then each quotient $\Lambda/I^i$ is discrete and $\Lambda=\lim_{i\geq 1}\Lambda/I^i$ as condensed rings so that \refsect{coefficients.limits.sect} applies. 
In the following, all limits are derived unless mentioned otherwise. 

An object $M\in \D(X,\Lambda)$ is called \emph{derived $I$-complete} if the natural map $M\r \lim (M\t_{\Lambda_X}(\Lambda/I^i)_X)$ is an equivalence. 
Let us temporarily denote by 
\[
\D_\lis(X,\widehat\Lambda)\subset \D_\cons(X,\widehat\Lambda) 
\]
the full subcategories of $\D(X,\Lambda)$ of objects $M$ which are derived $I$-complete and such that its reduction $M\t_{\Lambda_X}(\Lambda/I)_X$ is lisse, respectively constructible in $\D(X, \Lambda/I)$.
By \thref{discrete.comparison.prop} applied to the discrete ring $\Lambda/I$, the homotopy category of $\D_\cons(X,\widehat\Lambda)$ agrees with the category defined in \cite[Definition 6.5.1]{BhattScholze:ProEtale}. 

\lemm
\thlabel{reduction.I}
Let $M \in \D(X, \Lambda)$.
If $M \t_{\Lambda_X} (\Lambda_X / I)_X$ is lisse (respectively, constructible), then so is $M \t_{\Lambda_X} (\Lambda / I^i)_X$ for all $i \ge 1$.
\xlemm

\pf
Using \thref{discrete.comparison.prop} for $\Lambda/I$, one sees by induction on $i$ that each reduction $M\t_{\Lambda_X}(\Lambda/I^i)_X$ lies in the essential image of the functor $\nu^*\co \D(X_\et,\Lambda/I^i)\r \D(X_\proet,\Lambda/I^i)$ and is étale locally constant. 
Passing to a suitable étale covering of $X$, our claim for lisse sheaves reduces to the corresponding statement for $\Lambda/I^i$-modules and the nilpotent ideal $I/I^i$, see \StP{07LU}.  
For constructible sheaves $M$, it follows that any stratification witnessing the constructibility of $M\t_{\Lambda_X}(\Lambda/I)_X$ induces such a stratification for $M\t_{\Lambda_X}(\Lambda/I^i)_X$. 
\xpf

\prop
\thlabel{adic.comparison.prop}
There are equalities
\[
\D_\lis(X,\Lambda)= \D_\lis(X,\widehat \Lambda) \text{ and } \D_\cons(X,\Lambda)=\D_\cons(X,\widehat\Lambda).
\]  
as full subcategories of $\D(X,\Lambda)$.
\xprop

\pf
The canonical functor $\D_\lis(X,\widehat\Lambda)\r \lim \D_\lis(X,\Lambda/I^i)$ (and likewise for constructible sheaves) afforded by \thref{reduction.I} is fully faithful by \cite[Lemma 3.5.7 (2)]{BhattScholze:ProEtale}.
It is essentially surjective by \thref{limit.coefficients.lemm}, so we are done again using \thref{limit.coefficients.lemm}.
\xpf

\subsection{Classical approaches and some examples}
\label{sect--examples.comparison.l.adic}
Fix a prime $\ell$ and a finite extension $E / \Ql$. Its ring of integers $\calO_E$ has a unique maximal ideal $\mathfrak m$, for which we have $\calO_E = \lim \calO_E / \mathfrak m^n$.
The latter is a profinite topological ring, while $E = \calO_E[\ell^{-1}]$ carries the (usual) colimit topology so that $\calO_E\subset E$ is an open subring.

For a scheme $X$, let use denote by $\D_\ctf(X,\calO_E/\mathfrak m^n)$ the \ii-category of constructible étale $\calO_E/\mathfrak m^n$-sheaves of finite Tor-dimension. 
It is classical to consider the limit
\[
\Dbc(X,\calO_E)= \lim\D_\ctf(X,\calO_E/\mathfrak m^n),
\]
that is, the category of compatible systems of such objects. 
From here one usually passes to $E$-coefficients by inverting $\ell$:
\[
\Dbc(X,E)= \Dbc(X,\calO_E)\t_{\Perf_{\calO_E}}\Perf_{E} 
\]
This tensor product agrees with the idempotent completion of the localization $\Dbc(X,\calO_E)[\ell^{-1}]$. 
The following result gives the comparison with \cite[Exposé XIII, \S4]{Illusie:TravauxGabber}, \StP{0F4M} and \cite[Definition 6.8.8]{BhattScholze:ProEtale}:

\theo
\thlabel{classical.comparison.lemm}
For any qcqs scheme $X$, there are natural equivalences
\begin{align*}
\D_\cons(X,\calO_E/\mathfrak m^n)\cong & \ \D_\ctf(X,\calO_E/\mathfrak m^n), \\
\D_\cons(X,\calO_E)\cong & \ \Dbc(X,\calO_E) \text{ and } \\
\D_\cons(X,E)\cong & \ \Dbc(X,E).
\end{align*}
\xtheo
\pf
By \thref{discrete.comparison.prop} and \thref{ctf.comparison.rema}, the pullback of sheaves along $X_\proet\r X_\et$ induces equivalences $\D_\ctf(X,\calO_E/\mathfrak m^n)=\D_\cons(X_\et,\calO_E/\mathfrak m^n)\cong \D_\cons(X,\calO_E/\mathfrak m^n)$.
Passing to the limits, we obtain an equivalence $\Dbc(X,\calO_E)\cong \D_\cons(X,\calO_E)$ by \thref{limit.coefficients.lemm}.
Finally, we have a fully faithful embedding $\Dbc(X,E)\hr \D_\cons(X,E)$ by \thref{localization.lisse.cons.lemm} which is essentially surjective by \cite[Corollary 2.4]{HansenScholze:RelativePerversity}.
The latter uses \thref{insane.lattice.coro}. 
Alternatively, if $X$ is topologically Noetherian, one can use \thref{locally.constant.cor} to see that  $\D_\cons(X,E)$ agrees with the category defined in \cite[Definition 6.8.8]{BhattScholze:ProEtale}.
In this case, the result also follows from \cite[Proposition 6.8.14]{BhattScholze:ProEtale}.
\xpf

\rema
Likewise, one can show the following formula for any qcqs scheme $X$ with coefficients in the finite adeles $\bbA_{\bbQ,\on{f}}$:
\[
\Dcons(X,\bbA_{\bbQ,\on{f}})\cong \big(\mathrm{lim}\ \D_\ctf(X,\bbZ/n)\big)\t_{\Perf_{\widehat \bbZ}}\Perf_{\bbA_{\bbQ,\on{f}}}.
\]
\xrema

If $X$ is of finite type over a finite or separably closed field, then the homotopy category of $\D_\cons(X,\bbZ_\ell)$ is the $2$-limit of the categories $\D_\ctf(X,\bbZ/\ell^n)$ considered in \cite[Equation (1.1.2)]{Deligne:Weil2}, see also \cite[Section II.6]{KiehlWeissauer}.
For the relation with \cite{Ekedahl:Adic}, the reader is referred to \cite[Section 5.5]{BhattScholze:ProEtale}, in view of \refsect{adic.coefficients}. 
Also, $\D_\cons(X,\bbZ_\ell)$ is equivalent to the \ii-category defined in \cite[Definition 2.3.2.1]{GaitsgoryLurie:Weil} when $X$ is quasi-projective over an algebraically closed field.
We leave the details to the reader.

We finally turn to $\bar\bbQ_\ell$ or $\bar\bbZ_\ell$ coefficients.
Recall that $\bar\bbZ_\ell=\colim \calO_E$ and $\bar\bbQ_\ell=\colim E$ carry the colimit topology.
We recover the classical approach \cite[(1.1.3)]{Deligne:Weil2}, \cite[II.5, Appendix A]{KiehlWeissauer}:
\lemm
\thlabel{classical.comparison.colimit.lemm}
For any qcqs scheme $X$, there are natural equivalences
\begin{align*}
\colim_{E/\bbQ_\ell \;\on{finite}} \Dbc(X,\calO_E)\stackrel \cong \lr & \ \D_\cons(X,\bar\bbZ_\ell) \text{ and } \\
\colim_{E/\bbQ_\ell \;\on{finite}} \Dbc(X,E)\stackrel \cong \lr & \  \D_\cons(X,\bar\bbQ_\ell).
\end{align*}
\xlemm
\pf
This immediate from the above discussion using \thref{colimit.coefficients.lemm}.
\xpf

\rema[The $6$ functors]
\thlabel{six.functor.rema}
Under the usual finiteness, excellency and $\ell$-coprimality assumptions on the schemes, one obtains a $6$ functor formalism along the usual lines for the categories of constructible sheaves. 
In light of \thref{classical.comparison.lemm}, the reader is referred to the treatment in \cite[Section 6.7]{BhattScholze:ProEtale}.
\xrema

\section{Ind-lisse and ind-constructible sheaves}
\label{sect--ind.cons.sheaves}

Let $X$ be a qcqs scheme and $\Lambda$ a condensed ring. 
In this section, we impose the following finiteness assumption on the $\Lambda$-cohomological dimension:
there exists an integer $d_X\geq 0$ such that for all $p> d_X$,
all pro-\'etale affines $U=\lim U_i\in X_\proet$ and all sheaves $N\in \D(X,\Lambda)^\heartsuit$ of the form $N=\H^0(M)$ for some $M\in \D_\cons(X, \Lambda)$ we have the vanishing
\eqn
\label{cohomological.dimension.eq}
\H^p(U,N) \;=\; 0.
\xeqn
Examples include finite type schemes $X$ over finite or separably closed fields with coefficients $\Lambda$ either a discrete torsion ring, an algebraic extension $E/\bbQ_\ell$ or its ring of integers $\calO_E$, see \thref{Artin.vanishing}. 

Recall from \refsect{generalities.proetale} that the \ii-category $\D(X,\Lambda)$ admits small (co-)limits.  

\defi 
\thlabel{ind.lis.cons.defi}
A sheaf $M\in \D(X,\Lambda)$ is called \textit{ind-lisse}, respectively \textit{ind-constructible} if it is equivalent to a filtered colimit of lisse, respectively constructible $\Lambda$-sheaves.
\xdefi

The full subcategories of $\D(X,\Lambda)$ of ind-lisse, respectively ind-constructible sheaves are denoted by
\[
\D_\indlis (X,\Lambda)\subset \D_\indcons(X,\Lambda).
\]
Both \ii-categories are naturally commutative algebra objects in $\PrSt_{\Gamma(X,\Lambda)}$, that is, $\Gamma(X,\Lambda)$-linear symmetric monoidal stable presentable \ii-categories, see \thref{Gaitsgory.Lurie.comparison} for the properties stable and presentable. Recall the notion of compact objects \cite[Section 5.3.4]{Lurie:Higher}. 

\prop 
\thlabel{compact.objects}
An object $M\in \D_\indcons(X,\Lambda)$ is compact if and only if $M$ is constructible, and likewise for (ind-)lisse sheaves. 
Consequently, passing to compact objects induces equalities 
\[
\D_\indcons(X,\Lambda)^\omega=\D_\cons(X,\Lambda) \text{ and } \D_\indlis(X,\Lambda)^\omega=\D_\lis(X,\Lambda).
\]
\xprop

Before giving the proof, let us point out the following corollary which, for example, gives the comparison with the presentable categories of $\ell$-adic sheaves considered in \cite[\S 2.3.2]{GaitsgoryLurie:Weil}.
It is noteworthy since the condition in \eqref{cohomological.dimension.eq} is weaker than requiring constructible sheaves to be compact objects in $\D(X, \Lambda)$.

\coro
\thlabel{Gaitsgory.Lurie.comparison}
The inclusion $\D_\cons(X,\Lambda)\subset \D_\indcons(X,\Lambda)$ extends to a colimit-preserving equivalence
\[
\Ind\big(\D_\cons(X,\Lambda)\big) \stackrel{\cong} \lr \D_\indcons(X,\Lambda),
\]
and likewise for (ind-)lisse sheaves. 
\xcoro
\pf
By \cite[Proposition 5.3.5.11]{Lurie:Higher}, the functor is fully faithful because all objects of $\D_\cons(X,\Lambda)$ are compact in $\D_\indcons(X,\Lambda)$.
The essential surjectivity is immediate from the definition. 
\xpf

\rema
\thlabel{ind.cons.t.structure}
If $\Lambda$ is a discrete ring, then \thref{Gaitsgory.Lurie.comparison} together with \thref{discrete.comparison.prop} shows that $\D_\indcons(X,\Lambda)\cong \D(X_\et,\Lambda)$, see also \cite[Proposition 6.4.8]{BhattScholze:ProEtale}.
In particular, the natural t-structure on $\D(X,\Lambda)$ restricts to a t-structure on $\D_\indcons(X,\Lambda)$ even if $\Lambda$ is not t-admissible. 
\xrema

The proof of \thref{compact.objects} builds on the following lemma that crucially relies on assumption \eqref{cohomological.dimension.eq}, see \thref{cohomological.dimension.warn}:

\lemm
\thlabel{all.colimit.commutation.lemm}
The following functors commute with filtered colimits with terms in $\D_\cons(X,\Lambda)$:
\begin{enumerate}

\item 
\label{item--all.colimit.commutation.lemm.1}
$f_*\co \D(X,\Lambda)\r \D(Y,\Lambda)$ for map $f\co X\r Y$ between qcqs schemes satisfying \eqref{cohomological.dimension.eq};

\item 
\label{item--all.colimit.commutation.lemm.2}
$\underline\Hom_{\Lambda_X}(M,\str)\co \D(X,\Lambda)\r \D(X,\Lambda)$ for any qcqs $X$ scheme satisfying \eqref{cohomological.dimension.eq} and $M\in \D_\cons(X,\Lambda)$.
\end{enumerate}
In particular, under the conditions in \refit{all.colimit.commutation.lemm.2}, the functor 
\[
\Hom_{\Lambda_X}(M,\str)=\RG\big(X,\underline\Hom_{\Lambda_X}(M,\str)\big)\co \D(X,\Lambda)\r \Mod_{\Gamma(X,\Lambda)}
\]
commutes with such colimits as well. 
\xlemm

\pf
Assuming that part \refit{all.colimit.commutation.lemm.1} holds, the rest is proven analogously to \thref{colimit.commutation.lemm}.
For \refit{all.colimit.commutation.lemm.1}, let $\{M_i\}$ be a filtered system of constructible $\Lambda$-sheaves $M_i$ on $X$.
We need to show that the map 
\[
\colim \H^p\circ f_*(M_i)\r \H^p\circ f_*(\colim M_i), 
\]
is an equivalence in $\D(Y,\Lambda)^\heartsuit$ for every $p\in \bbZ$, $\H^p:=\tau^{\leq p}\circ \tau^{\geq p}$.
By \thref{colimit.commutation.lemm}, this holds true if all $M_i$ lie in $\D^{\geq n}(X,\Lambda)$ for some $n\in\bbZ$.
So it is enough to show that there exists $n\in \bbZ$ such that the map $\H^p\circ f_*(\tau^{\geq n}M)\r \H^p\circ f_*(M)$ is an isomorphism for all $M\in \D_\indcons(X,\Lambda)$.
Evaluating at any w-contractible pro-\'etale affine $V\in Y_\proet$, we obtain the map on pro\'etale cohomology groups
\[
\H^p\big(X\x_YV,\tau^{\geq n}M\big)\r \H^p\big(X\x_YV,M\big)
\]
Using the left-completeness of $\D(X,\Lambda)$ \cite[Proposition 3.3.3]{BhattScholze:ProEtale}, we are reduced to showing: there exists some integer $d_X\geq 0$ such that $\H^p(U,N)=0$ for all $p>d_X$, all qcqs $U\in X_\proet$ that admit an open cover by pro-\'etale affines and all $N\in \D(X,\Lambda)^\heartsuit$ of the form $N=\H^0(M)$ for some $M\in \D_\indcons(X,\Lambda)$. 
Since $\H^p(U,\str)$ commutes with filtered colimits in $\D(X,\Lambda)^\heartsuit$ (\thref{colimit.commutation.lemm}), we may assume that $N=\H^0(M)$ for some $M\in \D_\cons(X,\Lambda)$.
By an induction on the finite number pro-\'etale affines covering $U$, passing through the case of separated $U$ first, we may assume that $U=\lim U_i$ is pro-\'etale affine.
So the desired integer $d_X$ exists by our assumption \eqref{cohomological.dimension.eq}.
\xpf

\pf[Proof of \thref{compact.objects}]
We only treat (ind-)constructible sheaves as the argument for (ind-)lisse sheaves is completely analogous. 

To show that $M\in \D_\cons(X,\Lambda)$ is compact in $\D_\indcons(X,\Lambda)$, we need to show that the natural map
\[
\colim \Hom_{\Lambda_X}(M, N_j) \r \Hom_{\Lambda_X}(M,\colim N_j)
\]
is an equivalence for every filtered system $\{N_j\}$ of constructible $\Lambda$-sheaves. 
This follows from \thref{all.colimit.commutation.lemm}.

Conversely, let $M=\colim M_i\in \D_\indcons(X,\Lambda)$ be a compact object.
Then the identity $\id_M\co M\r M$ factors through some $M_i$, presenting $M$ as a direct summand of $M_i$.
As $\D_\cons(X,\Lambda)$ is idempotent complete, we see that $M$ is constructible.

\xpf

\lemm
\thlabel{Artin.vanishing}
The pair $(X,\Lambda)$ satisfies \eqref{cohomological.dimension.eq} in each of the following cases:
\begin{enumerate}
\item
\label{item--Artin.vanishing.1}
The scheme $X$ is of finite type over a finite or separably closed field and $\Lambda$ is either a discrete torsion ring, an algebraic field extension $E\supset \bbQ_\ell$ or its ring of integers $\calO_E$.
\item 
\label{item--Artin.vanishing.2}
The scheme $X$ is a qcqs scheme in characteristic $p>0$ and $\Lambda$ is either a discrete $p$-torsion ring, an algebraic field extension $E\supset \bbQ_p$ or its ring of integers $\calO_E$.
\end{enumerate}
\xlemm
\pf
Let $U=\lim U_i\in X_\proet$ be pro-\'etale affine, and let $N=\H^0(M)$ for some $M\in \D_\cons(X,\Lambda)$. 
The case where $\Lambda$ is an algebraic extension $E\supset \bbQ_\ell$ or its ring of integers $\calO_E$ are reduced to the case of finite extensions (\thref{classical.comparison.colimit.lemm}) and further to the case of finite discrete torsion rings (\thref{classical.comparison.lemm}).
It remains to treat the case where $\Lambda$ is a discrete torsion ring. 
Then $\H^p(U,N)=\colim \H^p(U_i,N)$ using \thref{discrete.comparison.prop} and the continuity of the \'etale site \StP{03Q4}. 
So part \refit{Artin.vanishing.1} follows from Artin vanishing (see \StP{0F0V}), noting that the abelian sheaf underlying $N$ is torsion and that separably closed fields (respectively, finite fields) have cohomological dimension $0$ (respectively, $1$).
For \refit{Artin.vanishing.2}, we claim that $\H^i(X,N)=0$ for all $i> 2$, affine schemes $X = \Spec R$ in characteristic $p>0$ and $p$-torsion abelian sheaves $N$. 
Using d\'evissage arguments \StPd{09Z4}{03SA}, it suffices to consider a constructible $\bbF_p$-sheaf $N$.
By topological invariance of the \'etale site, we may assume $R$ to be perfect. 
There is a fully faithful t-exact functor from the \ii-category $\Dcons(X, \Fp)$ to the derived category of modules over $R[F]$ defined in \cite[Notation~2.1.5]{BhattLurie:Riemann}, see Theorem 12.1.5 there. 
This functor sends the constant sheaf $\Fp$ to $R$, which has a length two resolution by projective $R[F]$-modules, see \cite[Section 3]{BhattLurie:Riemann}. 
This shows the claim.
\xpf

\coro 
\thlabel{descent.ind.cons}
Let $U_\bullet \r X$ be a hypercover such that $ U_n$ is quasi-compact \'etale for all $n\geq 0$.
Then the natural functor
\[
\D_{\indcons}(X,\Lambda) \stackrel\cong\lr \Tot\big(\D_{\indcons}(U_\bullet,\Lambda)\big)
\]
is an equivalence.
If $U_n\r X$ is finite \'etale, then the same holds for ind-lisse sheaves. 
\xcoro
\pf 
By the descent equivalence $\Tot\big(\D(U_\bullet, \Lambda)\big) = \D(X, \Lambda)$, there are full inclusions
\[
\D_{\indcons}(X,\Lambda) \;\subset\; \Tot\big(\D_{\indcons}(U_\bullet,\Lambda)\big) \;\subset\; \D(X,\Lambda).
\]
It remains to show: if $M\in \D(X,\Lambda)$ such that $M|_U$ is ind-constructible, then $M$ is so. 
We denote $j_\bullet\co U_\bullet\r X$.
For every $M\in \D(X,\Lambda)$, we have a canonical equivalence $|(j_\bullet)_{!}\circ  j^*_\bullet M| \xrightarrow{\sim} M$. 
Since $U\r X$ is finitely presented \'etale, the same is true for every $j_{n}\colon U_n \rightarrow X$. 
In particular, each functor $(j_n)_!$ preserves the constructible subcategories (\thref{preservation.constructibility}). 
So, if $j_{n}^{*}M$ can be written as a filtered colimit of constructible objects, then the same is true for $(j_{n})_!j_{n}^{*}M$. 
This shows $M\in \D_{\indcons}(X,\Lambda)$.

Finally, if $U\r X$ is finite \'etale, then each functor $(j_n)_!$ preserves the lisse subcategories (\thref{preservation.constructibility}), and we can proceed as before.
\xpf

\begin{appendix}

\section{The condensed shape of the pro\'{e}tale topos}
\label{sect--condensed.shape}

Let $X$ be a qcqs scheme, and $\Lambda$ a condensed ring. 
In this appendix, we explain how the formalism developed above gives a simple realization of lisse $\Lambda$-sheaves on $X$ as representations of the condensed shape associated to $X_\proet$ valued in perfect $\Lambda$-modules. 
This is related to the stratified shape developed in \cite{Barwick.Glasman.Haine.exodromy, BarwickHaine:Pyknotic, wolf2020proetale}. 
This appendix is not used throughout the manuscript. 
We include the material as it gives another point of view on the \ii-categories of lisse sheaves $\D_\lis(X,\Lambda)$ introduced in \refsect{generalities.proetale}.

We take the following definition of the condensed shape which is similar to the classical definitions of Artin--Mazur--Friedlander.
Denote by $\mathrm{HC}(X)$ the \ii-category of hypercoverings in $X_{\proet}$ whose objects consist of hypercovers $U_{\bullet}\to X$ in $X$ with $U_{n}$ qcqs for all $n\geq 0$. 
For precise definitions, the reader is referred to \cite{dugger.hollander.isaksen.hypercovers} and \cite[\textsection{5}]{hoyois.higher.galois.theory}. 
The \ii-category $\mathrm{HC}(X)$ is cofiltered by \cite[Proposition 5.1]{dugger.hollander.isaksen.hypercovers}. 
The \textit{condensed shape of $X_\proet$} is the condensed animated set
\[
\Picond(X) \defined \lim_{U_\bullet\in \mathrm{HC}(X)} |\pi_0(U_\bullet)|,
\]
where the geometric realization and the limit are taken in the \ii-category $\Cond(\Ani)$ of condensed anima (also called spaces, Kan complexes, or \ii-groupoids). 
Here we identify a profinite set with the associated condensed set under the Yoneda embedding. 

Let $\mathrm{HC}^{w}(X) \subset \mathrm{HC}(X)$ denote the full subcategory consisting of hypercovers $U_\bullet$ so that, for every $n\geq 0$, the scheme $U_n$ is w-contractible. 
Since every hypercover can be refined by one consisting of w-contractible schemes, the inclusion $\mathrm{HC}^{w}(X) \subset \mathrm{HC}(X)$ is co-initial. 
In particular, the natural map
\[
\Picond(X) \stackrel\cong\lr \lim_{U_\bullet \in \mathrm{HC}^{w}(X)} |\pi_0(U_\bullet)|.
\]
is an equivalence.
Since covers of w-contractible qcqs objects split, all the condensed sets in the colimit are actually equivalent. 
As the category $\mathrm{HC}^{w}$ is cofiltered, we get an equivalence of condensed anima
\eqn\label{shape.w-contractible.cover}
\Picond(X) \stackrel\cong\lr |\pi_0(U_\bullet)|,
\xeqn
for every $U_\bullet\in \mathrm{HC}^{w}(X)$.

For the condensed ring $\Lambda$, we define a condensed \ii-category $\Perf^{\mathrm{cond}}_{\Lambda}$ by the assignment
\begin{align*}
    \Perf^{\mathrm{cond}}_{\Lambda}\colon \{\text{extremally disconnected profinite sets}\}^{op} &\r \Cat_\infty^\perf\\
    S & \mapsto \Perf_{\Gamma(S,\Lambda)}.
\end{align*}
Note that $\Perf^{\mathrm{cond}}_{\Lambda}$, when extended to a hypersheaf on $*_\proet$, is simply the hypersheaf $S \mapsto \D_\lis(S,\Lambda)$, see \thref{dualizable.objects.contractible.lem} and \thref{lisse.constructible.hyperdescent.coro}.

The category of condensed  \ii-categories has a canonical enrichment over $\Cat_\infty$. 
Namely, for every \ii-category $K$, we have the associated $\underline{K}\in \Cond(\Cat_\infty)$ as the presheaf on the category of extremally disconnected profinite sets
\[
S=\lim S_i \mapsto \colim_{i} \Fun(S_i,K).
\]
Then, for condensed \ii-categories $\mathcal{C}$, $\mathcal{D}$, the mapping object $\Fun^{\cts}(\mathcal{C},\mathcal{D})\in \Cat_{\infty}$ is characterized by the existence of natural equivalences
\[
\Hom_{\Cat_\infty}\big(K,\Fun^{\cts}(\mathcal{C},\mathcal{D})\big) \cong \Hom_{\Cond(\Cat_\infty)}\big(\underline{K}\times \mathcal{C},\mathcal{D}\big)
\]
with $K\in \Cat_\infty$. 
In order to simplify notation, we identify every profinite set $S$ with the associated discrete condensed \ii-category. 
Then, for any extremally disconnected profinite set $S$, we have the equivalence
\[
\Perf_{\Gamma(S,\Lambda)} \cong \Fun^{\cts}(S,\Perf^{\mathrm{cond}}_{\Lambda}),
\]
using the Yoneda embedding.

\prop
\thlabel{lisse.sheaves.via.condensed.shape}
Let $X$ be a qcqs scheme and $\Lambda$ a condensed ring. 
Then there is a canonical equivalence
\eqn\label{profinite.galois.equivalence}
\D_\lis(X,\Lambda) \cong \Fun^{\cts}\big(\Pi_{\mathrm{cond}}(X),\Perf_{\Lambda}^{\mathrm{cond}} \big).
\xeqn
\xprop
\pf 
Let $U_\bullet\in \mathrm{HC}^{w}(X)$. 
By descent (\thref{lisse.constructible.hyperdescent.coro}), we have an equivalence
\eqn \label{descent.proof.of.condensed.shape}
\D_\lis(X,\Lambda) \cong \Tot \left(\Perf_{\Gamma(U_\bullet,\Lambda)}\right) \cong \Tot \left(\Perf_{\Gamma(\pi_0(U_\bullet), \Lambda)}\right),
\xeqn
using $\Gamma(U_n,\Lambda)=\Gamma(\pi_0(U_n), \Lambda)$ for $n\geq 0$.
Since each $U_n$ is w-contractible qcqs, each profinite set $\pi_0(U_n)$ is extremally disconnected. 
We get an equivalence
\[
\Perf_{\Gamma(\pi_0(U_\bullet),\Lambda)} \cong \Fun^{\cts}\big(\pi_0(U_\bullet),\Perf_{\Lambda}^{\mathrm{cond}} \big).
\]
Then \eqref{descent.proof.of.condensed.shape} becomes 
\[
\D_\lis(X,\Lambda) \cong \tot\left(\Fun^{\cts}\big(\pi_0(U_\bullet),\Perf_{\Lambda}^{\mathrm{cond}} \big)\right) \cong 
\Fun^{\cts}\big(|\pi_0(U_\bullet)|,\Perf_{\Lambda}^{\mathrm{cond}} \big).
\]
This construction is functorial in $U_\bullet$ and gives a canonical equivalence 
\[
\D_\lis(X,\Lambda) \cong \Fun^{\cts}\big(\Pi_{\mathrm{cond}}(X),\Perf_{\Lambda}^{\mathrm{cond}} \big).
\]
\xpf 

This proposition allows us to draw direct connections between functors from the condensed shape and more familiar versions when one restricts the possible class of rings.
First, if $\Lambda$ is finite and discrete, then profinite shape theory in the sense of \cite[Appendix~E]{Lurie:SAG} gives an equivalence:
\[
\D_\lis(X,\Lambda) \cong \Fun(\Pi_{\proet}(X),\Perf_\Lambda),
\]
where $\Pi_{\proet}(X)$ denotes the profinite shape of the topos $X_\proet$. Namely, the profinite completion of the shape. Here we use $\Fun$ as in \cite{hoyois.higher.galois.theory} to denote the category of functors between pro-categories. 
Now consider a condensed ring $\Lambda$ associated to a Noetherian ring complete with respect to the adic topology for some ideal $I\subset \Lambda$ with $\Lambda/I$ finite. 
Then \thref{limit.coefficients.lemm} implies
\[
\D_\lis(X,\Lambda) \cong \lim_{n} \Fun(\Pi_{\proet}(X),\Perf_{\Lambda/I^n}) \cong \Fun(\Pi_{\proet}(X),\Perf_{\Lambda}),
\]
where on the right hand side both arguments are considered as pro-categories. 

\end{appendix}

\bibliographystyle{alphaurl}
\bibliography{bib}

\end{document}